\documentclass[11pt, a4paper, twoside]{article}
\usepackage[a4paper]{anysize}

\usepackage{amssymb}
\usepackage{booktabs}
\usepackage{amsmath}
\usepackage{amsfonts, amsthm, graphicx, stmaryrd, fancyhdr}
\usepackage{enumerate}

\newcommand{\ud}{\mathrm{d}}
\newcommand{\E}{\mathbb{E}}
\newcommand{\F}{\mathcal{F}}
\newcommand{\I}{\int_{t}^T}
\newcommand{\Lp}{\mathbb{L}^p}
\newcommand{\Lzwei}{\mathbb{L}^2}
\newcommand{\N}{\mathbb{N}^*}

\newcommand{\Q}{\mathbb{Q}}
\newcommand{\R}{\mathbb{R}}

\newcommand{\om}{{\omega}}

\newcommand{\cE}{\mathcal{E}}
\newcommand{\cF}{\mathcal{F}}

\newcommand{\cL}{\mathcal{L}}

\newcommand{\cO}{\mathcal{O}}

\newcommand{\cR}{\mathcal{R}}
\newcommand{\cS}{\mathcal{S}}

\newcommand{\ebs}{e^{\beta s}}
\newcommand{\ebt}{e^{\beta t}}
\newcommand{\ebT}{e^{\beta T}}
\newcommand{\Ft}{\mathcal{F}_t}

\newcommand{\be}{\begin{eqnarray*}}
\newcommand{\ee}{\end{eqnarray*}}
\newcommand{\ben}{\begin{eqnarray}}
\newcommand{\een}{\end{eqnarray}}

\theoremstyle{plain}
\newtheorem{theo}{Theorem}[section]
\newtheorem{lemma}[theo]{Lemma}

\newtheorem{corollary}[theo]{Corollary}

\theoremstyle{definition}
\newtheorem{defi}[theo]{Definition}

\title{Classical and Variational Differentiability of BSDEs with Quadratic Growth}
\author{Stefan Ankirchner and Peter Imkeller and Gon\c calo dos Reis
\\ Institut f\"ur Mathematik\\ Humboldt-Universit\"at zu Berlin\\
Unter den Linden 6\\ 10099 Berlin\\ Germany }
 
\begin{document}

\maketitle

\begin{abstract}
We consider Backward Stochastic Differential Equations (BSDEs) with
generators that grow quadratically in the control variable. In a
more abstract setting, we first allow both the terminal condition
and the generator to depend on a vector parameter $x$. We give
sufficient conditions for the solution pair of the BSDE to be
differentiable in $x$. These results can be applied to systems of
forward-backward SDE. If the terminal condition of the BSDE is given
by a sufficiently smooth function of the terminal value of a forward
SDE, then its solution pair is differentiable with respect to the
initial vector of the forward equation. Finally we prove sufficient
conditions for solutions of quadratic BSDEs to be differentiable in
the variational sense (Malliavin differentiable).
\end{abstract}



%

\section*{Introduction}

Problems of stochastic control treated by the crucial tool of
\emph{backward stochastic differential equations (BSDEs)} have been
encountered in many areas of application of mathematics in recent
years. A particularly important area is focused around optimal
hedging problems for contingent claims in models of financial
markets. Recently, a special class of hedging problems in incomplete
financial markets has been considered in the area where finance and
insurance concepts meet. At this interface problems of
\emph{securitization} arise, i.e. insurance risk is transferred to
capital markets. One particularly interesting risk source is given
by climate or environmental hazards affecting insurance companies or
big branches of the economy that depend on weather such as
agriculture and fishing, transportation and tourism. The public
awareness of climate hazards such as floods or hurricanes is
continually increasing with the intensity of the discussion about
irreversible changes due to human impact.\par\smallskip

BSDEs typically appear in the following setting. On a financial
market some small investors are subject to an external risk source
described for instance by weather or climate influences. There may
also be big investors such as re-insurance companies that depend
in a possibly different way on the same risk source. In this
situation market incompleteness stems from the external risk not
hedgeable by the market assets. One may complete the market either
by making the external risk tradable through the introduction of
an insurance asset traded among small agents, or by introducing a
risk bond issued by a big agent. In this setting, treating the
utility maximization problem for the agents under an equilibrium
condition describing basically market clearing for the additional
assets, leads to the determination of the market price of external
risk through a BSDE which in case of exponential utility turns out
to be quadratic in the control variable (see \cite{06HM},
\cite{05CHIM} and \cite{04CIM}). Alternatively, instead of
maximizing utility with respect to exponential utility functions
we might minimize risk measured by the \emph{entropic risk
measure}. In this setting we again encounter a BSDE with quadratic
nonlinearity, of the type

$$Y_t = \xi + \int_t^T f(s,Y_s,Z_s) ds - \int_t^T Z_s d W_s,\quad
0\le t\le T,
$$

where $W$ is a finite-dimensional Wiener process of the same
dimension as the control process $Z$, with a generator $f$ that
depends at most quadratically on $Z$, and a bounded terminal
condition $\xi.$ In the meantime, the big number of papers published
on general BSDEs is rivalled by the number of papers on BSDEs of this
type of nonlinearity. For a more complete list of references see
\cite{05CSTV} or \cite{00Kob}. In particular, there are papers in
which the boundedness condition on $\xi$ is relaxed to an
exponential integrability assumption, or where the stochastic
integral process of $Z$ is supposed to be a BMO
martingale.\par\medskip

In a particularly interesting case the terminal variable $\xi$ is
given by a function $g(X^x_T)$ at terminal time $T$ of the solution
process $X$ of a forward SDE
$$X^x_t = x + \int_0^t b(s, X^x_s) ds + \int_0^t \sigma(s, X^x_s) d
W_s,\quad 0\le t\le T,$$ with initial vector $x\in \mathbb{R}$.
Similarly, the driver $f$ may depend on the diffusion dynamics of
$X^x$. Via the famous link given by the generalized Feynman-Kac
formula, systems as the above of forward-backward stochastic
differential equations are seen to yield a stochastic access to
solve nonlinear PDE in the viscosity sense, see \cite{00Kob}.\par\smallskip

In this context, questions related to the regularity of the
solutions $(X^x, Y^x, Z^x)$ of the stochastic forward-backward
system in the classical sense with respect to the initial vector $x$ or in the
sense of the stochastic calculus of variations (Malliavin
calculus) are frequently encountered. Equally, from a more
analytic point of view also questions of smoothness of the
viscosity solutions of the PDE associated via the Feynman-Kac link
are seen to be very relevant.\par\smallskip

For instance, Horst and M\"uller (see \cite{06HM}) ask for
existence, uniqueness and regularity of a global classical
solution of our PDE from the \emph{analytic point of view}. Not
attempting a systematic approach of the problem, they use the
natural access of the problem by asking for smoothness of the
solutions of the stochastic system in terms of the stochastic
calculus of variations. But subsequently they work under the
restrictive condition that the solutions of the BSDE have bounded
variational derivatives, which is guaranteed only under very
restrictive assumptions on the coefficients.\par\smallskip

The question of smoothness of the \emph{stochastic solutions} in
the parameter $x$ arises for instance in an approach of cross
hedging of environmental risks in \cite{05AIP}. Here the setting is
roughly the one of an incomplete market generated by a number of
big and small agents subject to an external (e.g. climate related)
risk source, and able to invest in a given capital market. The
risk exposure of different types of agents may be negatively
correlated, so that typically one type profits from the risky
event, while at the same time the other type suffers. Therefore
the concept of hedging one type's risk by transferring it to the
agents of the other type in a \emph{cross hedging} context makes
sense. Mathematically, in the same way as described above, it
leads to a BSDE of the quadratic type, the solution $(Y^x, Z^x)$
of which depends on the initial vector $x$ of a forward equation
with solution $X^x.$ Under certain assumptions, the cross-hedging
strategy can be explicitly given in a formula depending crucially
on $x$, and in which the sensitivity with respect to $x$ describes
interesting quality properties of the strategy.\par\medskip

In this paper, we tackle regularity properties of the solutions
$(Y^x,Z^x)$ of BSDEs of the quadratic type such as the two previously
sketched in a systematic and thorough way. Firstly, the particular
dependence on the starting vector $x$ of the forward component of a
forward-backward system will be generalized to the setting of a
terminal condition $\xi(x)$ depending in a smooth way to be
specified on some vector $x$ in a certain Euclidean state space. We
both consider the smoothness with respect to $x$ in the classical
sense, as well as the smoothness in the sense of Malliavin's
calculus.\par\smallskip

The common pattern of reasoning in order to tackle smoothness
properties of any kind starts with a priori estimates for difference
and differential quotients, or for infinite dimensional gradients in
the sense of variational calculus. In the estimates, these
quantities are related to corresponding difference and differential
quotients or Malliavin gradients of the terminal variable and the
driver. To obtain the a priori estimates, we make use to changes of
probability of the Girsanov type, by which essentially nonlinear
parts of the driver are eliminated. Since terminal conditions in our
treatment are usually bounded, the exponential densities in these
measure changes are related to $BMO$ martingales. Known results
about the inverse H\"older inequality allow to show that as a
consequence the exponential densities are $r$-integrable for some
$r>1$ related to the $BMO$ norm. This way we are able to reduce
integrability properties for the quantities to be estimated to a
natural level. In a second step, the a priori inequalities are used
to derive the desired smoothness properties from corresponding
properties of driver and terminal condition. To the best of our
knowledge, only Malliavin differentiability results of this type
have been obtained so far, with strong conditions on the
coefficients restricting generality considerably (see
\cite{06HM}).\par\smallskip

The paper is organized as follows. In section \ref{sec:prelim} we
fix the notation and recall some process properties needed in the
proofs of the main body of the paper. Section \ref{diffsection}
contains the main results on classical differentiability. In
sections \ref{secesti}, \ref{secpriori} and
\ref{section:priori.bsde} we give a priori bounds for classes of
non-linear BSDEs. Section \ref{sec:proof:diff} contains the proofs of
the theorems stated in Section \ref{diffsection}. Section
\ref{sec:forward-backward} is devoted to the application of the
proven results to the forward-backward SDE setting. In Section
\ref{sec:malliavin} we state and prove the Malliavin
differentiability results.

\section{Preliminaries}\label{sec:prelim}
Throughout this paper let $(\Omega, \F, P)$ be a complete
probability space and $W=(W_t)_{t\geq 0}$ a $d-$dimensional Brownian
motion. Let $\{\Ft\}_{t\geq 0}$ denote the natural filtration
generated by $W$, augmented by the $P-$null sets of $\F$.

Let $T > 0$, $\xi$ be an $\cF_T$-measurable random variable and
$f:\Omega \times [0,T] \times \R \times \R^d \to \R$. We will
consider Backward Stochastic Differential Equations (BSDEs) of the
form \ben\label{bsde0} Y_t = \xi + \int_t^T  f(t,Y_t,Z_t)\ud t -
\int_t^T Z_t\ud W_t. \een As usual we will call $\xi$ the {\em
terminal condition} and the function $f$ the {\em generator} of
the BSDE (\ref{bsde0}). A {\em solution} consists of a pair
$(Y,Z)$ of adapted processes such that (\ref{bsde0}) is satisfied.
To be correct we should write $\int_t^T \langle Z_t,\ud
W_t\rangle$ or $\sum_{i=1}^d \int_t^T Z_s^i dW_s^i$ instead of
$\int_t^T Z_t\ud W_t$, since $W$ and $Z$ are $d-$dimensional
vectors; but for simplicity we use this notation as it is without
ambiguity. It is important to know which process spaces the
solution of a BSDE belongs to. We therefore introduce the
following notation for the spaces we will frequently use. Let
$p\in [1, \infty]$. Then, for $m\in \N$
\begin{itemize}
\item $\Lp(\R^m)$ is the space of all progressively measurable
processes $(X_t)_{t\in[0,T]}$ with values in $\R^m$ such that
$\|X_t\|_{\Lp}^p = \E[\left( \int_0^T |X_s|^2 \ud s
\right)^{p/2} ]<\infty.$ 
\item  $\cR^p(\R^m)$ is the space
of all measurable processes $(X_t)_{t\in[0,T]}$ with values in
$\R^m$ such that $\| X \|_{\cR^p}^p = \E[\left( \sup_{t \in [0,T]}
|X_t| \right)^{p}]< \infty$. Note that $\cR^\infty(\R^m)$ is the
space of bounded measurable processes.

\item $H^p(\R^m)$ is the class of all local martingales $X$ such
that $\lVert X\lVert_{H^p}^p=\E^{P}[\langle
X\rangle_T^{\frac{p}{2}}]<\infty$.

\item $L^p(\R^m; P)$ is the space of $\cF_T$-measurable random
variables $X:\Omega\mapsto\R^m$ such that $\lVert
X\lVert^p_{L^p}=\E^P[|X|^p]<\infty$. We will omit reference to the
space or the measure when there is no ambiguity.

\end{itemize}
Furthermore, we use the notation $\partial_t=\frac{\partial}{\partial t}$,
$\nabla =(\frac{\partial}{\partial x_1},\cdots,\frac{\partial}{\partial x_d})$ for $(t,x)\in[0,T]\times\R^d$.

Suppose that the generator satisfies, for $a\geq 0$ and $b,c>0$
\ben\label{quadraticgene}
|f(t,x,y,z)|
\leq a(1+b|y|)+\frac{c}{2}|z|^2.
\een
Kobylanski has shown in \cite{00Kob} that if $\xi$ is bounded and
the generator $f$ satisfies (\ref{quadraticgene}), then there exists
a solution  $(Y,Z)\in \cR^\infty\times \mathbb{L}^2$. Moreover, it
follows from the results in \cite{06Mor}, that in this case the
process $Z$ is such that the stochastic integral process relative to
the Brownian motion $\int_0^\cdot Z dW$ is a so-called Bounded Mean
Oscillation (BMO) martingale.

Since the BMO property is crucial for the proofs we present in this
paper we recall its definition and some of its basic properties. For
an overview on BMO martingales see \cite{kazamaki}.

\begin{defi}[BMO]
Let $M$ be a uniformly integrable $(\F_t)$-martingale satisfying
$M_0=0$. For $1\leq p < \infty$ set
\[
\lVert M\lVert_{BMO_p}=\sup_{\tau \textrm{ stopping time}} \Big(\E\Big[|M_{\infty}-M_\tau|^p|\F_\tau\Big]\Big)^{1/p}.
\]
The normed linear space $\{M:\lVert M\lVert_{BMO_p}<\infty\}$ with
norm $\lVert M\lVert_{BMO_p}$ is denoted by $BMO_p$. If we want to
stress the measure $P$ we are referring to we will write BMO($P$).

It can be shown that for any $p$, $q \in [1,\infty]$ we have
$\mathrm{BMO}_p = \mathrm{BMO}_q$ (see \cite{kazamaki}). Therefore
we will often omit the index and simply write BMO for the set of BMO
martingales.
\end{defi}

In the following Lemma we state the properties of BMO martingales we
will frequently use.
\begin{lemma}[Properties of BMO martingales]\label{bmoeigen}$\phantom{121}$
\begin{itemize}
\item[1)] Given a BMO martingale $M$ with quadratic variation $\langle M\rangle$, its stochastic exponential
\[\cE(M)_T = \exp\{M_T-\frac{1}{2}\langle M\rangle_T\}\]
has integral $1$, and thus the measure defined by $dQ = \cE(M)_T dP$
is a probability measure.

\item[2)] Let $M$ be a BMO martingale relative to the measure $P$.
Then the process $\hat{M}= M - \langle M\rangle$ is a BMO martingale
relative to the measure $Q$ (see Theorem 3.3 in \cite{kazamaki}).

\item[3)] For any BMO Martingale, it is always possible to find a $p>1$ such that $\cE(M) \in L^p$, i.e.
if $\lVert M\lVert_{BMO_2}<\Psi(p)$, then $\cE(M) \in L^p$ (see for
example Theorem 3.1 \cite{kazamaki}). Where $\Psi(x)=\Big\{1+\frac{1}{x^2}\log\frac{2x-1}{2(x-1)}\Big\}^{\frac{1}{2}}-1$ for all $1<x<\infty$
and verifies $\lim_{x\to1^+}
\Psi(x)=\infty$ and $\lim_{x\to\infty} \Psi(x)=0$.

\end{itemize}
\end{lemma}
\section{Differentiability of quadratic BSDEs in the classical sense} \label{diffsection}
Suppose that the terminal condition and the generator of a quadratic
BSDE depend on the Euclidean parameter set $\R^n$ for some $n \in
\N$. We will show that the smoothness of the terminal condition and
the generator is transferred to the solution of the BSDE
\begin{equation}\label{parabsde}
Y_t^{x}=\xi(x)-\int_t^T Z_s^{x}\ud W_s+\int_t^T
f(s,x,Y_s^x,Z_s^{x})\ud s,\quad x\in\R^n,
\end{equation}
where terminal condition and generator are subject to the following
conditions
\begin{enumerate}[({C}1)]
\item\label{para1} $f:\Omega \times [0,T]\times \R^n
\times\R\times\R^d\to\R$ is an adapted measurable function such
that $f(\omega,t,x,y,z)=l(\omega,t,x,y,z)+\alpha |z|^2$, where
$l(\omega, t, x, y, z)$ is globally Lipschitz in $(y,z)$ and
continuously differentiable in $(x,y,z)$; for all $r\ge1$ and
$(t,y,z)$ the mapping $\R^d \to L^r$, $x \mapsto l(\om, t,x,y,z)$
is differentiable and for all $x\in\R^n$ \be
\lim_{x'\to x} \E^P\Big[\Big( \int_0^T |l(s, x', Y_s^x, Z_s^x)-l(s, x, Y_s^x, Z_s^x)|\ud s \Big)^r\Big]= 0 \qquad  \textrm{ and }\\
\lim_{x'\to x} \E^P\Big[\Big( \int_0^T |\frac{\partial}{\partial x} l(s, x', Y_s^{x'}, Z_s^{x'})- \frac{\partial}{\partial x}l(s, x, Y_s^{x}, Z_s^{x})|\ud s \Big)^r\Big]= 0,
\ee
\item\label{para2} the random variables $\xi(x)$ are $\F_T-$adapted and for every compact
set $K \subset \R^n$ there exists a constant $c\in \R$ such that
$\sup_{x \in K} \| \xi(x) \|_\infty \le c$; for all $p\ge1$ the
mapping $\R^n \to L^p$, $x \mapsto \xi(x)$ is differentiable with
derivative $\nabla \xi$.
\end{enumerate}
If (C\ref{para1}) and (C\ref{para2}) are satisfied, then there
exists a unique solution $(Y^x,Z^x)$ of Equation (\ref{parabsde}).
This follows from Theorems 2.3 and 2.6 in \cite{00Kob}. We will
establish two differentiability results for the pair $(Y^x,Z^x)$
in the variable $x$. We first consider differentiability of the
vector valued map $$x \mapsto (Y^x, Z^{x})$$ with respect to the Banach space topology
defined on $\cR^p(\R^1) \times \Lp(\R^d)$. This will be stated in Theorem \ref{paradiffbsde2}.
A slightly more stringent result will be obtained in the subsequent Theorem \ref{paradiffbsde}. Here, we
consider pathwise differentiability of the maps
$$x \mapsto (Y^x_t(\omega), Z^{x}_t(\omega))$$ in the usual sense, for almost all pairs $(\omega, t)$. In both cases, the
derivatives will be identified with $(\nabla Y^x,\nabla Z^x)$ solving the BSDE
\begin{equation}\label{pararealdiff}
\begin{array}{lll}
\nabla Y_t^{x}
&=& \nabla \xi(x)-\int_t^T \nabla Z_s^{x}\ud W_s\\
&&  +\int_t^T \left[ \partial_x l(s,x,Y_s^{x},Z_s^{x}) + \partial_y l(s,x,Y_s^{x},Z_s^{x}) \nabla Y_s^{x} +
\partial_z l(s,x,Y_s^{x},Z_s^{x})\nabla Z_s^{x} + 2 \alpha Z^x_s \nabla Z_s^{x}\right] \ud
s.
\end{array}
\end{equation}
We emphasize at this place that it is not immediate that this BSDE
possesses a solution. In fact, without considering it as a
component of a system of BSDEs also containing the original
quadratic one, it can only be seen as a linear BSDE with global,
but random (and not bounded) Lipschitz constants.
\begin{theo}\label{paradiffbsde2}
Assume (C\ref{para1}) and (C\ref{para2}). Then for all $p \ge 1$, the function $\R^n \to \cR^p(\R^1) \times \Lp(\R^d)$,
$x \mapsto (Y^x, Z^{x})$, is differentiable, and the derivative is a solution of the BSDE (\ref{pararealdiff}).
\end{theo}
Under slightly stronger conditions one can show the existence of a
modification of $Y^x$ which is $P$-a.s. differentiable as a mapping
from $\R^n$ to $\R$. Let $e_i=(0, \ldots, 1, \ldots, 0)$ be the unit
vector in $\R^n$ where the $i$th component is $1$ and all the other
components $0$. For $x\in\R^n$ and $h \not=0$ let $\zeta(x,h,e_i) =
\frac1h [\xi(x+h e_i) - \xi(x)]$. For the existence of
differentiable modifications we will assume that
\begin{enumerate}
\item[(C3)]
for all $p \ge 1$ there exists a constant $C >0$ such that for all
$i \in\{1, \ldots, n\}$, $x$, $x' \in \R^n$ and $h, h' \in
\R\setminus \{0\}$
\[ \E\Big[|\xi(x+he_i) - \xi(x'+h'e_i)|^{2p} + |\zeta(x,h,e_i) - \zeta(x',h',e_i)|^{2p}\Big] \le C(|x-x'|^2 +
|h-h'|^2)^p.\]

\end{enumerate}
\begin{theo}\label{paradiffbsde}
Suppose, in addition to the assumptions of Theorem
\ref{paradiffbsde2}, that (C3) is satisfied and that $l(t,x,y,z)$
and its derivatives are globally Lipschitz continuous in $(x,y,z)$.
Then there exists a function $\Omega \times [0, T] \times \R^n \to
\R^{1+d}$, $(\om, t, x) \mapsto (Y^{x}_t, Z^x_t)(\om)$, such that
for almost all $\om$, $Y^x_t$ is continuous in $t$ and continuously
differentiable in $x$, and for all $x$, $(Y^x_t, Z^x_t)$ is a
solution of (\ref{parabsde}).
\end{theo}

\section{Moment estimates for linear BSDEs with stochastic Lipschitz generators} \label{secesti}

By formally deriving a quadratic BSDE with generator satisfying
(C\ref{para1}) and (C\ref{para2}) we obtain a linear BSDE with a
stochastic Lipschitz continuous generator. The Lipschitz constant
depends on the second component of the solution of the original
BSDE. In order to show differentiability, we start deriving a priori
estimates for this type of linear BSDE with stochastic Lipschitz
continuous generator. For this purpose, we first need to show that
the moments of the solution can be effectively controlled. Therefore
this section is devoted to moment estimates of solutions of BSDEs of
the form
\begin{equation} \label{bsde-lin-001}
U_t = \zeta - \int_t^T V_s dW_s + \int_t^T \left[ l(s,U_s, V_s)+H_s
V_s + A_s\right] \ud s.
\end{equation}
We will make the following assumptions concerning the drivers:
\begin{itemize}
\item[(A1)] For all $p \ge 1$, $\zeta$ is $\F_T-$adapted and we have $\zeta \in L^p(\R^1)$,
\item[(A2)] $H$ is a predictable $\R^d-$valued process, integrable with respect to $W$, such that $\int H  dW$ is a
BMO-martingale,
\item[(A3)] $l:\Omega \times [0,T] \times \R \times \R^d \to \R$ is such that for all $(u,v)$,
the process $l(\om,t,u,v)$ is $(\cF_t)$-predictable and there exists a constant $M > 0$ such that for all $(\om,t,u,v)$,
\[ |l(\omega, t, u,v)|\leq M(|u| +|v|), \]
\item[(A4)] $A$ is a measurable adapted process such that for all $p \ge 1$ we
have $\E[\big(\int_0^T |A_s|\ud s\big)^{p}] < \infty$.
\end{itemize}
Moreover, we assume that $(U,V)$ is a solution of (\ref{bsde-lin-001}) satisfying
\begin{itemize}
\item[(A5)] $[\int_0^T  U_s^2 |V_s|^2 \ud s]^{\frac{1}{2}}$ and $\int_0^T  |U_s A_s| \ud s$ are $p$-integrable for all $p \ge 1$.

\end{itemize}
Under the assumptions (A1), (A2), (A3), (A4) and (A5) one obtains
the following estimates.
\begin{theo}[Moment estimates]\label{apriori.estimate}
Assume that (A1)-(A5) are satisfied. Let $p > 1$ and $r>1$ such that
$\mathcal{E}(\int H dW)_T \in L^r(P)$. Then there exists a constant
$C > 0$, depending only on $p$, $T$, $M$ and the BMO-norm of $\int H
dW)$, such that with the conjugate exponent $q$ of $r$ we have
\ben\label{moment.estimate.theo.01} \E^{P} \Big[\sup_{t\in[0,T]}
|U_t|^{2p} \Big]+\E^{P} \Big[\Big(\int_0^T |V_s|^2\ud s \Big)^p
\Big]
 &\leq&
    C\E^{P}\Big[\,| \zeta|^{2pq^2} + \Big(\int_0^T |A_s| \ud s\Big)^{2pq^2}\Big]^{\frac{1}{q^2}}.
\een
Moreover we have
\ben\label{moment.estimate.theo.02}
\E^{P} [\int_0^T |U_s|^{2}\ud s ] + \E^P[\int_0^T |V_s|^2\ud s]&\leq&
    C\E^{P}\Big[\,| \zeta|^{2q^2} + \Big(\int_0^T |A_s|^2 \ud s\Big)^{q^2}
    \Big]^{\frac{1}{q^2}}.
\een
\end{theo}
The proof is divided into several steps. First let $\beta >0$ and
observe that by applying It\^o's formula to $e^{\beta t} U_t^2$ we
obtain \be e^{\beta t} U_t^2&=&e^{\beta T} U_T^2-2\I e^{\beta s} U_s
V_s \ud W_s \\&&+ \I e^{\beta s}\big[-\beta U_s^2
+2U_s\big(l(s,U_s,V_s)+H_sV_s + A_s\big)-|V_s|^2\big]\ud s. \ee By
(A2), the auxiliary measure defined by $ Q = \mathcal{E}(H \cdot
W)_T \cdot P$ is in fact a probability measure. Then $\hat W_t = W_t
- \int_0^t H_s \ud s$ is a $Q$-Brownian motion, and \be e^{\beta t}
U_t^2 &\leq& e^{\beta T} U_T^2-2\I e^{\beta s} U_s V_s \ud \hat{W}_s
\\&&+\I e^{\beta s}\big[(-\beta+2M) U_s^2 +2M|U_s||V_s|-|V_s|^2 + |U_s
A_s|\big]\ud s\nonumber \ee By choosing $\beta=M^2+2M$, we obtain
\ben\label{ineq.base} e^{\beta t} U_t^2 +\I e^{\beta s}
(M|U_s|-|V_s|)^2\ud s &\leq&e^{\beta T} U_T^2-2\I e^{\beta s} U_s
V_s \ud \hat{W}_s + \I e^{\beta s} |U_s A_s| \ud s.
\een We therefore first prove moment estimates under the measure
$Q$.
\begin{lemma}\label{apriori.estimate.lemma}
For all $p > 1$ there exists a constant $C$, depending only on $p$, $T$ and $M$, such that
\ben\label{moment.estimate.lem.02}
\E^{Q} \Big[\sup_{t\in[0,T]} |U_t|^{2p} \Big]+ \E^{Q} \left[ \left(\int_0^T |V_s|^2\ud s \right)^p \right]
&\leq&  C\E^{Q}\Big[\,| \zeta|^{2p}+\Big( \int_0^T  |A_s| \ud s \Big)^{2p}\Big].
\een Moreover we have
\ben\label{moment.estimate.lem.01}
\E^{Q} \Big[\int_0^T |U_s|^{2}\ud s \Big] + \E^{Q} \left[\int_0^T |V_s|^2\ud s \right]&\leq&
    C\E^{Q}\Big[\,| \zeta|^{2} + \int_0^T  |A_s|^2 \ud s \Big].
\een
\end{lemma}

\begin{proof}
Throughout this proof let $C_1, C_2, \ldots$, be constants depending only on $p$, $T$ and $M$.

Inequality (\ref{ineq.base}) implies \ben\label{ineq.base1} e^{\beta
t} U_t^2 &\leq&e^{\beta T} U_T^2-2\I e^{\beta s} U_s V_s \ud {\hat
W}_s + \I e^{\beta s} |U_s A_s| \ud s, \een and (A5) together with
the existence of the $r$th moment for $\mathcal{E}(\int H dW)_T$
yield $\int_0^T U^2_s |V_s|^2 \ud s \in L^1(Q)$. Hence, since
$e^{\beta t} U_t^2$ is $(\F_t)$-adapted, \ben\label{151206-1}
e^{\beta t} U_t^2&\leq& e^{\beta T}\E^Q\Big[ |\zeta|^2 + \I e^{\beta
s} |U_s A_s| \ud s| \F_t\Big]. \een Integrating both sides and using
Young's inequality, we obtain \be \E^Q[\int_0^T  U_s^2\ud s]&\leq&
C_1 \E^Q[ |\zeta|^2 + \int_0^T  |U_s A_s| \ud s ]
\\&\leq&C_1 \E^Q[ \zeta^2 + 2C_1\int_0^T  |A_s|^2 \ud s ]+\frac{1}{2} \E^Q[\int_0^T U_s^2\ud s],
\ee and hence \ben\label{151206-2} \E^Q[\int_0^T  U_s^2\ud s]\leq
C_2 \E^Q[ |\zeta|^2 +\int_0^T  |A_s|^2 \ud s ]. \een Inequality
(\ref{151206-1}), (A5) and Doob's $L^p$ inequality imply for $p > 1$
\be
     \E^Q[\sup_{t\in[0,T]} | U_t|^{2p}]
     &\leq& C_3 \E^{Q}\Big[ \Big(|\zeta|^{2} + \int_0^T |U_s A_s| \ud s\Big)^p \Big]\\
     &\leq& C_4 \E^{Q}\Big[ |\zeta|^{2p} + \Big(\sup_{t\in[0,T]} |U_t| \int_0^T |A_s| \ud s\Big)^p \Big].
\ee
By Young's inequality, $(\sup_{t\in[0,T]} |U_t|^p) (\int_0^T |A_s| \ud s)^p \leq \frac{1}{2C_4}\sup_{t\in[0,T]} |U_t|^{2p} +2C_4(\int_0^T |A_s| \ud s)^{2p}$, and hence
\ben\label{firstpart}
     \E^Q[\sup_{t\in[0,T]} | U_t|^{2p}]
     \leq C_5 \E^{Q}\Big[ |\zeta|^{2p} + \Big(\int_0^T |A_s| \ud s\Big)^{2p} \Big].
\een
In order to complete the proof, note that (\ref{ineq.base}) implies
\ben\label{aux.piori.01}
&&\I e^{\beta s} |V_s|^2\ud s\nonumber \\
&&\qquad\leq e^{\beta T} U_T^2-2\I e^{\beta s} U_s V_s \ud \hat{W}_s
+2\I \ebs M|U_s||V_s|\ud s+\I \ebs |U_s||A_s|\ud s.
\een
By Young's inequality, $2 \I e^{\beta s} M|U_s||V_s|\ud s \le \frac12 \I e^{\beta s} |V_s|^2\ud s +
8 M^2 \I e^{\beta s} U_s^2 \ud s$, and hence
\be
\frac{1}{2}\E^Q\Big[\int_0^T e^{\beta s} |V_s|^2\ud s\Big]
 &\leq&\E^Q[e^{\beta T} U_T^2+8M^2\int_0^T \ebs U_s^2\ud s+ \int_0^T \ebs U_s^2+ \ebs|A_s|^2\ud s \Big]\\
  &\leq& C_6\E^{Q}[ \zeta^2 +\int_0^T |A_s|^2\ud s]
\ee
which, combined with (\ref{151206-2}) leads to the desired Inequality (\ref{moment.estimate.lem.01}).

Equation (\ref{aux.piori.01}), Young's inequality, Doob's $L^p$-inequality and the Burkholder-Davis-Gundy
inequality imply
\be
&&\E^{Q} \left[ \left(\int_0^T e^{\beta s}|V_s|^2\ud s \right)^p \right] \\
&\leq& C_7 \E^{Q} \Big[|\zeta|^{2p}+\Big(T \sup_{t\in[0,T]}{e^{\beta t}U_t^2} \Big)^p
+\Big(\int_0^T e^{\beta s}U_sV_s\ud \hat{W}_s \Big)^p + \Big( \int_0^T \ebs |U_s||A_s|\ud s \Big)^p\Big]\nonumber\\
&\leq&
C_8 \E^{Q} \Big[|\zeta|^{2p}+\sup_{t\in[0,T]}{e^{\beta tp}|U_t|^{2p}}
 + \Big(\int_0^T e^{2\beta s} U_s^2|V_s|^2\ud s\Big)^{\frac{p}{2}}
       +\sup_{t\in[0,T]}|U_t|^{p} \Big( \int_0^T \ebs|A_s|\ud s\Big)^p      \Big]\nonumber\\
&\leq&
C_8 \E^{Q} \Big[|\zeta|^{2p}+\sup_{t\in[0,T]}{e^{\beta tp}|U_t|^{2p}}
 + \Big( \sup_{t\in[0,T]}\ebt U_t^2\Big)^{\frac{p}{2}} \Big(\int_0^T \ebs |V_s|^2\ud s\Big)^{\frac{p}{2}}\\
&& \qquad \qquad+ \sup_{t\in[0,T]}|U_t|^{2p} + \Big( \int_0^T \ebs |A_s|\ud s\Big)^{2p} \Big]
\ee
By Young's inequality,
\[\Big( \sup_{t\in[0,T]}\ebt U_t^2\Big)^{\frac{p}{2}} \Big(\int_0^T \ebs |V_s|^2\ud s\Big)^{\frac{p}{2}} \le
2C_8 \Big(\sup_{t\in[0,T]}\ebt U_t^2\Big)^{p} + \frac{1}{2C_8}\Big(\int_0^T \ebs |V_s|^2\ud s\Big)^{p},\]
which implies
\be
\E^{Q} \left[ \left(\int_0^T e^{\beta s}|V_s|^2\ud s \right)^p \right]
&\leq& C_9 \E^{Q} \Big[|\zeta|^{2p}+\sup_{t\in[0,T]}|U_t|^{2p}+\Big( \int_0^T |A_s|\ud s\Big)^{2p}\Big]\\
&\leq& C_{10} \E^Q\Big[|\zeta|^{2p}+\Big( \int_0^T |A_s|\ud s\Big)^{2p}\Big].
\ee
Thus, with Inequality (\ref{firstpart}), the proof is complete.
\end{proof}
\begin{proof}[Proof of Theorem \ref{apriori.estimate}]
Notice that by the second statement of Lemma \ref{bmoeigen}, the process $\int H d\hat{W}=\int H dW-\int_0^\cdot H_s^2 \ud s$
belongs to BMO($Q$), and hence $-\int H d\hat{W}$ also. Moreover, $\cE(\int H dW)^{-1} = \cE(-\int H d\hat{W})$.
Consequently, by the third statement of Lemma \ref{bmoeigen}, there exists an $r > 1$ such that $\cE(H \cdot W)_T \in L^r(P)$
and $\cE(H \cdot W)^{-1}_T \in L^r(Q)$. Throughout let $D = \max\{ \| \cE(H \cdot W)_T \|_{L^r(P)}, \|
\cE(H \cdot W)^{-1}_T\|_{L^r(Q)} \}$. H\"older's inequality and Lemma \ref{apriori.estimate.lemma} imply that for the conjugate exponent
$q$ of $r$
we have
\be
\E^P[\sup_{s\in[0,T]} |U_s|^{2p}]&=& \E^Q[\cE(H \cdot W)^{-1}_T \sup_{s\in[0,T]} |U_s|^{2p}]\
\leq\ D \E^Q[ \sup_{s\in[0,T]} |U_s|^{2pq}]^{\frac1q} \\
&\leq& C_1\,D\,\E^{Q}\Big[ |\zeta|^{2pq} + \Big(\int_0^T |A_s| \ud s\Big)^{2pq} \Big]^{\frac1q}\\
&=&C_1\,D\, \E^P[\cE(H \cdot W)_T \Big( |\zeta|^{2pq} + \Big(\int_0^T |A_s| \ud s\Big)^{2pq}\Big)]^{\frac1q}\\
&\leq& C_2 \,D^{\frac{1+q}{q}} \E^P[|\zeta|^{2pq^2} + \Big(\int_0^T |A_s| \ud s\Big)^{2pq^2}]^{\frac{1}{q^2}},
\ee
where $C_1, C_2$ represent constants depending on $p, M, T$ and the $BMO$ norm of $\int H dW$.
Similarly, with another constant $C_3$,
$\E^P[\int_0^T |V_s|^{2p} \ud s] \leq C_3\,D^{\frac{1+q}{q}} \E^P[|\zeta|^{2pq^2} + \Big(\int_0^T |A_s|
\ud s\Big)^{2pq^2}\Big)]^{\frac{1}{q^2}}$, and hence (\ref{moment.estimate.theo.01}).
By applying the same arguments to (\ref{moment.estimate.lem.01}) we finally get (\ref{moment.estimate.theo.02}).
\end{proof}
%
%
%
%
%
\section{A priori estimates for linear BSDEs with stochastic Lipschitz constants}\label{secpriori}
In this section we shall derive a priori estimates for the variation
of the linear BSDEs that play the role of good candidates for the
derivatives of our original BSDE. These will be used to prove
continuous differentiability of the smoothly parametrized solution
in subsequent sections. Let $(\zeta, H, l_1, A)$ and $(\zeta', H',
l_2, A')$ be parameters satisfying the properties (A1), (A2), (A3)
and (A4) of Section \ref{secesti} and suppose that $l_1$ and $l_2$
are globally Lipschitz continuous and differentiable in $(u,v)$. Let
$(U,V)$ resp. $(U',V')$ be solutions of the linear BSDE
\begin{equation} \label{bsde.compar}
U_t = \zeta - \int_t^T V_s dW_s + \int_t^T [l_1(s,U_s, V_s)+H_s V_s+A_s] \ud s
\end{equation}
resp.
\be U'_t = \zeta' - \int_t^T V'_s dW_s + \int_t^T [l_2(s,U'_s,
V'_s)+H'_s V'_s+A'_s] \ud s \ee both satisfying property (A5).
Throughout let $\delta U_t = U_t - U'_t$, $\delta V_t = V_t - V'_t$,
$\delta \zeta = \zeta - \zeta'$, $\delta A_t= A_t-A'_t$ and $\delta
l(t,u,v) =l_1(t,u,v)-l_2(t,u,v)$.
\begin{theo}[A priori estimates]\label{comparison.theo}
Suppose we have for all $\beta \ge 1$,  $\int_0^T \delta U_s^2 |\delta V_s|^2 \ud s \in L^\beta(P)$ and
$\int_0^T |\delta U_s \delta A_s| \ud s \in L^\beta(P)$. Let $p \ge 1$ and $r>1$
such that $\mathcal{E}(\int H' dW)_T \in L^r(P)$. Then there exists a constant $C > 0$, depending
only on $p$, $T$, $M$ and the BMO-norm of $\int H' dW$, such that with the conjugate exponent $q$ of $r$ we have
\be
&&\!\!\!\E^P\Big[\sup_{t\in[0,T]} |\delta U_t|^{2p}\Big]+\E^P\Big[\Big(\int_0^T |\delta V_s|^2\ud s\Big)^p\Big]\\
&&\qquad\leq
C \Big\{ \E^P\Big[|\delta \zeta|^{2pq^2}+\Big(\int_0^T  |\delta l(s,U'_s, V'_s)+\delta A_s|\ud s\Big)^{2pq^2}\Big]^{\frac{1}{q^2}}  \\
 & &  \qquad\qquad + \ (\E^P[|\zeta|^{2pq^2}+\big(\int_0^T |A_s|\ud s \big)^{2pq^2}])^\frac{1}{2q^2} \E^P\Big[\Big(\int_0^T |H_s-H'_s|^2\ud\nonumber s\Big)^{2pq^2}\Big]^{\frac{1}{2q^2}}\Big\}
\ee
\end{theo}
We proceed in the same spirit as in the preceding section.
Before proving Theorem \ref{comparison.theo} we will show a priori estimates with respect to the
auxiliary probability measure $Q$ defined by $Q=\mathcal{E}(\int H' dW)_T \cdot P$.
Note that $\hat W_t = W_t - \int_0^t H'_s \ud s$ is a $Q$-Brownian motion.
\begin{lemma}\label{comparison.ineqs.Qmeasure}
Let $p>1$. There exists a constant $C > 0$, depending only on $p$, $T$ and $M$, such that

\ben\label{ineq.lemma.01}
\E^Q\Big[\sup_{t \in[0,T]}|\delta U_t|^{2p}\Big]
&\leq&C\Big\{\E^Q\Big[ |\delta \zeta|^{2p}+\Big(\int_0^T  |\delta l(s,U'_s, V'_s)+\delta A_s|\ud s\Big)^{2p}\Big]\\
   &&+ \left(\E^Q\Big[|\zeta|^{2p}+\big(\int_0^T |A_s|\ud s \big)^{2p}\Big]\right)^\frac12\nonumber
   \E^Q\Big[\Big(\int_0^T |H_s-H'_s|^2\ud s\Big)^{2p}\Big]^{\frac{1}{2}}\Big\},
\een
\ben\label{ineq.lemma.02}
\E^Q[\Big(\int_0^T |\delta V_s|^2\ud s\Big)^p]
&\leq&
C \Big\{ \E^Q\Big[|\delta \zeta|^{2p}+\Big(\int_0^T |\delta l(s,U'_s, V'_s)+\delta A_s|\ud s\Big)^{2p}\Big]\\
&&+ \left(\E^Q\Big[|\zeta|^{2p}+\big(\int_0^T |A_s|\ud s \big)^{2p}\Big]\right)^\frac12
\E^Q\Big[\Big(\int_0^T |H_s-H'_s|^2\ud s\Big)^{2p}\Big]^{\frac{1}{2}}\Big\}.\nonumber
\een
\end{lemma}

\begin{proof}
The difference $\delta U$ satisfies
\be
\delta U_t &=& \delta \zeta - \int_t^T \delta V_s \ud W_s
     + \int_t^T [(H_s V_s - H'_s V'_s) +l_1(s,U_s,V_s)-l_2(s,U'_s,V'_s)+\delta A_s]\ud s\\
 &=& \delta \zeta - \int_t^T \delta V_s \ud W_s+\I [l_1(s,U'_s,V'_s)-l_2(s,U'_s,V'_s)+ H'_s\delta V_s+\delta A_s]\ud s\\
&&     +\int_t^T [(H_s-H'_s)V_s +l_1(s,U_s,V_s)-l_1(s,U'_s,V'_s)]\ud
s. \ee
Let $\beta >0$. Applying It\^o's formula to $e^{\beta t} \delta
U_t^2, t\ge 0$, yields the equation
\ben\label{eq.aux.005} e^{\beta
t} \delta U_t^2&=&e^{\beta T} \delta U_T^2-2\I e^{\beta s} \delta
U_s \delta V_s \ud W_s
+2\I e^{\beta s} \delta U_s\,H'_s \,\delta V_s\ud s \nonumber\\
&&+\I e^{\beta s}\Big[-\beta \delta U_s^2 - |\delta V_s|^2
+2\big(l_1(s,U_s,V_s)-l_1(s,U'_s,V'_s)\big)\delta U_s\Big]\ud s\nonumber\\
&&+2\I e^{\beta s} \delta U_s (H_s-H'_s)V_s\ud s+2\I\ebs \delta U_s (\delta l_s +\delta A_s)\ud s,
\een
where $\delta l_s=l_1(s,U'_s,V'_s)-l_2(s,U'_s,V'_s)$.
Using the Lipschitz property of $l_1$ we obtain
\be
e^{\beta t} \delta U_t^2&\leq&\ebT \delta U_T^2
+\I e^{\beta s}\Big[(-\beta+2M)\delta U_s^2 - |\delta V_s|^2 +2M|\delta U_s|\,|\delta V_s|\Big]\ud s\\
&&+2\I e^{\beta s} \delta U_s [(H_s-H'_s)V_s +\delta l_s + \delta A_s]\ud s -2\I e^{\beta s} \delta U_s \delta V_s \ud \hat{W}_s.
\ee
If $\beta=(M^2+2M)$, then
\ben
\label{eq.aux.002}
e^{\beta t} \delta U_t^2+\I e^{\beta s}(M|\delta U_s|-|\delta V_s|)^2\ud s
&\leq&\ebT \delta U_T^2+2\I e^{\beta s} \delta U_s [(H_s-H'_s)V_s + \delta l_s + \delta A_s]\ud s\nonumber\\
&&-2\I e^{\beta s} \delta U_s \delta V_s \ud \hat{W}_s. \een

We will now derive the desired estimates from Equation
(\ref{eq.aux.002}). First observe that by taking conditional
expectations, we get \be e^{\beta t} \delta U_t^2 &\leq&\ebT
\E^Q\left[ \delta U_T^2+2\I e^{\beta s} \delta U_s [(H_s-H'_s)V_s +
\delta l_s + \delta A_s]\ud s \big| \F_t \right]. \ee Let $p> 1$. Then
for some constants $C_1, C_2, \ldots$, depending on $p$, $T$ and
$M$, we obtain \be \sup_{t\in[0,T]}|\delta U_t|^{2p} &\leq&\! C_1
\sup_{t\in[0,T]}\left\{ \Big(\E^Q\big[ |\delta U_T|^{2}|\F_t\big]
+\E\big[\int_0^T |\delta U_s [(H_s-H'_s)V_s + \delta l_s + \delta
A_s]| \ud s\big|\F_t\big]\Big)^{p}\right\} \ee and by Doob's $L^p$
inequality we get \be \E^Q[\sup_{t\in[0,T]}|\delta U_t|^{2p}]
&\leq&C_2\left\{ \E^Q\Big[ |\delta U_T|^{2p}] +\E[\Big(\int_0^T
|\delta U_s [(H_s-H'_s)V_s + \delta l_s + \delta A_s]| \ud
s\Big)^{p}]\right\}. \ee By using Young's and H\"older's
inequalities we have \ben\label{eq.aux.010}
&&\E^Q\Big[\Big(\int_0^T |\delta U_s [(H_s-H'_s)V_s + \delta l_s + \delta A_s]|\ud s\Big)^{p}\Big]\nonumber\\
&&\leq C_3 \E^Q\left\{\sup_{t\in[0,T]}|\delta U_t|^{p}
\left[ \Big(\int_0^T |H_s-H'_s|^2\ud s\Big)^{\frac{p}{2}}\Big( \int_0^T |V_s|^2\ud s\Big)^{\frac{p}{2}} + \Big( \int_0^T |\delta l_s + \delta A_s| \ud s\Big)^p  \right] \right\} \nonumber\\
&&\leq \frac{1}{2C_4}\E^Q\Big[\sup_{t\in[0,T]} |\delta
U_t|^{2p}\Big]  \nonumber\\
&&\qquad +4C_4\E^Q\Big[\Big(\int_0^T |H_s-H'_s|^2\ud
s\Big)^{p}
\Big(\int_0^T |V_s|^2\ud s)\Big)^{p} + \Big( \int_0^T |\delta l_s + \delta A_s| \ud s\Big)^{2p}\Big]\nonumber\\
&&\leq
\frac{1}{2C_4}\E^Q\Big[\sup_{t\in[0,T]} |\delta U_t|^{2p} \Big]+ C_5 \Big\{ \E^Q\Big( \int_0^T |\delta l_s + \delta A_s| \ud s\Big)^{2p}\nonumber\\
&&\qquad\qquad+\E^Q\Big[\Big(\int_0^T (H_s-H'_s)^2\ud
s\Big)^{2p}\Big]^{\frac{1}{2}} \E^Q\Big[\Big(\int_0^T  |V_s|^2\ud
s\Big)^{2p}\Big]^{\frac{1}{2}} \Big\}. \een
Therefore, we may further estimate
\be
\E^Q[\sup_{t\in[0,T]} |\delta U_t|^{2p}]
&\leq&C_6 \Big\{ \E^Q[\, |\delta \zeta|^{2p}]+\E^Q[\big(\int_0^T  |\delta l_s + \delta A_s|\ud s\big)^{2p}] \\
&&+\E^Q\Big[\Big(\int_0^T |H_s-H'_s|^2\ud s\Big)^{2p}\Big]^{\frac{1}{2}}
\E^Q\Big[\Big(\int_0^T |V_s|^2\ud s\Big)^{2p}\Big]^{\frac{1}{2}} \Big\}.
\ee
Due to Lemma \ref{apriori.estimate.lemma}, $\E^Q\Big[\Big(\int_0^T |V_s|^2\ud s\Big)^{2p}\Big]^{\frac{1}{2}}
\leq C_7 \E^Q\Big[|\zeta|^{2p}+(\int_0^T |A_s|\ud s)^{2p}\Big]^{\frac{1}{2}}<\infty$, which implies the $\delta U_s$ part of Inequality (\ref{ineq.lemma.01}).

In order to prove the second inequality, note that (\ref{eq.aux.002}) also implies
\ben\label{eq.aux.007}
\I e^{\beta s}|\delta V_s|^2\ud s
&\leq&\ebT \delta U_T^2+2\I e^{\beta s} \delta U_s [(H_s-H'_s)V_s + \delta l_s + \delta A_s]\ud s\nonumber\\
&&+2M\I e^{\beta s}|\delta U_s|\,|\delta V_s|\ud s
-2\I e^{\beta s} \delta U_s \delta V_s \ud \hat{W}_s.
\een

Equation (\ref{eq.aux.007}), Doob's $L^p$-inequality and the Burkholder-Davis-Gundy
inequality imply 
\be
\E^Q[\Big(\int_0^T |\delta V_s|^2\ud s\Big)^p]
&\leq&
C_8 \bigg\{(\E^Q\Big[|\delta \zeta|^{2p}+\int_0^T |\delta U_s|^{2p}\ud s\Big]+\E^Q[\Big(\int_0^T  \delta U_s^2 \delta |V_s|^2 \ud s\Big)^{\frac{p}{2}}] \\
&&+\E^Q[\big(\int_0^T |\delta U_s [(H_s-H'_s)V_s + \delta l_s + \delta A_s]|\ud s\big)^p] \bigg\}.
\ee
Consequently, Young's inequality allows to deduce
\be
\E^Q[\Big(\int_0^T |\delta V_s|^2\ud s\Big)^p]
&\leq&
C_9 \bigg\{ \E^Q\Big[|\delta \zeta|^{2p}+\int_0^T |\delta U_s|^{2p}\ud s\Big]+\E^Q[\sup_{t\in[0,T]}  |\delta U_t|^{2p}]
\\ &&+\E^Q[\big(\int_0^T |\delta U_s [(H_s-H'_s)V_s + \delta l_s + \delta A_s]|\ud s\big)^p] \bigg\}.
\ee
Finally, (\ref{ineq.lemma.01}) and (\ref{eq.aux.010}) imply
\be
\E^Q[\Big(\int_0^T |\delta V_s|^2\ud s\Big)^p]
&\leq&
C_{10} \E^Q\Big[|\delta \zeta|^{2p}+\Big(\int_0^T |\delta l_s + \delta A_s|\ud s\Big)^{2p}\Big]\\
&&\qquad+C_{10}\E^Q[|\zeta|^{2p}+(\int_0^T |A_s|\ud s)^{2p}]^{\frac{1}{2}}\E^Q\Big[\Big(\int_0^T |H_s-H'_s|^2\ud s\Big)^{2p}\Big]^{\frac{1}{2}}
\ee
and hence the proof is complete.
\end{proof}
\begin{proof}[Proof of Theorem \ref{comparison.theo}]
This can be deduced from Lemma \ref{comparison.ineqs.Qmeasure} with arguments similar to those of
Theorem \ref{apriori.estimate}. We just have to invoke Lemma \ref{bmoeigen}.
\end{proof}

\section{A priori estimates for quadratic BSDEs}\label{section:priori.bsde}

Consider the two quadratic BSDEs
\begin{equation} \label{bsde-lin}
Y_t = \xi - \int_t^T Z_s dW_s + \int_t^T [l_1(s,Y_s,Z_s)+\alpha Z_s^2] \ud s
\end{equation}
and
\begin{equation}
Y'_t = \xi' - \int_t^T Z'_s dW_s + \int_t^T [l_2(s,Y'_s,Z'_s)+\alpha (Z'_{s})^2] \ud s,
\end{equation}
where $\xi$ and $\xi'$ are two bounded $\cF_T$-measurable random
variables, and $l_1$ and $l_2$ are globally Lipschitz and
differentiable in $(y,z)$. Put now $\delta Y_t=Y_t-Y'_{t}$,
$\delta Z_t=Z_{t}-Z'_{t}$, $\delta \xi=\xi - \xi'$ and $\delta
l=l_1-l_2$. The a priori estimates we shall prove next will serve
for establishing (moment) smoothness of the solution of the
quadratic BSDE with respect to a parameter on which the terminal
variable depends smoothly. Note first that by boundedness of $\xi$
and $\xi'$ we have that both $\int Z dW$ and $\int Z' dW$ are
$BMO$ martingales, so that we may again invoke the key Lemma
\ref{bmoeigen}.

\begin{theo} \label{apriori-quadra}
Suppose that for all $\beta \ge1$ we have $\int_0^T |\delta
l(s,Y_s,Z_s)| \ud s \in L^\beta(P)$. Let $p >1$ and choose $r>1$
such that $\mathcal{E}(\alpha (Z_s+ Z'_s)\cdot W)_T \in
L^r(P)$. Then there exists a constant $C > 0$, depending only on
$p$, $T$, $M$ and the BMO-norm of $(\alpha \int(Z_s+ Z'_s)
dW)$, such that with the conjugate exponent $q$ of $r$ we have  \be
&&\E^P \Big[\sup_{t\in[0,T]} |\delta Y_t|^{2p} \Big]+\E^{P} \left[
\left(\int_0^T |\delta Z_s|^2\ud s \right)^{p} \right]\\&&
\qquad \qquad \leq C \,
\left( \E^{P}\Big[\,|\delta \xi|^{2p q^2} +(\int_0^T |\delta
l(s,Y_s, Z_s)|\ud s)^{2pq^2} \Big]\right)^{\frac{1}{q^2}}. \ee
Moreover we have \be \E^P[ \int_0^T |\delta Y_s|^2\ud s]+\E^{P}
\big[ \int_0^T |\delta Z_s|^2\ud s \big]
 &\leq &C \left( \E^{P}\left[\,|\delta \xi|^{2q^2} +\Big(\int_0^T |\delta l(s,Y_s, Z_s)|\ud s\Big)^{2q^2}\right] \right)^{\frac{1}{q^2}}.
\ee
\end{theo}
We give only a sketch of the proof since the arguments are very
similar to the ones used in the proofs in Sections \ref{secesti} and
\ref{secpriori}.

First observe that \be \delta Y_t &=& \delta \xi - \int_t^T \delta Z_s
\ud W_s + \int_t^T [l_1(s,Y_s,Z_s)-l_1(s,Y'_s,Z'_s)+\delta
l(s,Y'_s,Z'_s) + \alpha(Z_s + Z_s') \delta Z_s] \ud s. \ee
By applying It\^o's formula to $e^{\beta t}|\delta Y_t|^2$ we obtain 
\begin{eqnarray} \label{ito}
&& e^{\beta t}|\delta Y_t|^2 - e^{\beta T}|\delta Y_T|^2 \nonumber\\
&&= 2\int_t^T e^{\beta s}\delta Y_s\Big( l_1(s,Y_s,Z_s)-l_1(s,Y'_s,Z'_s) +\delta l(s,Y'_s,Z'_s)\Big) \ud s
   -2\I e^{\beta s}(\beta|\delta Y_s|^2+|\delta Z_s|^2)\ud s  \nonumber \\
&&\qquad+ 2\int_t^T e^{\beta s}\delta Y_s \alpha (Z_s+Z'_s)\delta Z_s \ud s
 - 2\I e^{\beta s}\delta Y_s\delta Z_s\ud W_s.
\end{eqnarray}
We start with a priori estimates under the auxiliary probability
measure $Q$ defined  by $Q = \mathcal{E}(\alpha \int
(Z_s+Z'_s) dW) \cdot P$. Note that $\tilde W_t = W_t - \int_0^t
\alpha (Z_s + Z'_s) ds$ is a $Q$-Brownian motion.

Let $\beta>0$. Equality (\ref{ito}) and
the Lipschitz property of $l_1$  yield
\begin{eqnarray*}
e^{\beta t}|\delta Y_t|^2 &\leq& e^{\beta T} |\delta \xi|^2
-2\I e^{\beta s}\,\delta Y_s\,\delta Z_s\ud \tilde W_s+2\int_t^T \ebs \delta Y_s \delta l(s,Y'_s,Z'_s)\ud s\\
&&+\I e^{\beta s}\Big( (-\beta+2M)|\delta Y_s|^2-|\delta
Z_s|^2+2M|\delta Y_s||\delta Z_s|\Big)\ud s.
\end{eqnarray*}
By choosing $\beta=M^2+2M$ we obtain the general inequality
\ben\label{eqBMO001}
&&e^{\beta t}|\delta Y_t|^2+\I e^{\beta s}(|\delta Z_s|-M|\delta Y_s|)^2\ud s\nonumber\\
&&\qquad\qquad\leq e^{\beta T}|\delta \xi|^2 +\int_t^T \ebs \delta
Y_s \delta l(s,Y'_s,Z'_s)\ud s -2\I e^{\beta s} \delta Y_s \, \delta
Z_s \ud \tilde W_s. \een
Note that the process $\int_0^t e^{\beta
s}\,\delta Y_s\, \delta Z_s\ud \tilde W_s$ is a strict
martingale because $\delta Y_s$ is bounded and $(\delta Z \cdot
\tilde W)$ is BMO relative to $Q$.

Notice that Equation (\ref{eqBMO001}) is of similar but simpler form
than Equation (\ref{eq.aux.002}). This is because the $(H_s-H'_s)$
term in (\ref{eq.aux.002}) has been completely absorbed by the
Girsanov measure change. As a consequence, following the proof of Lemma \ref{comparison.ineqs.Qmeasure}, we obtain
the following estimates:
\begin{lemma}
For all $p >1$ there exists a constant $C >0$, depending only on $p$, $M$ and $T$, such that
\be
\E^Q \Big[\sup_{t\in[0,T]} |\delta Y_t|^{2p} \Big]+\E^{Q} \left[ \left(\int_0^T |\delta Z_s|^2\ud s \right)^{p} \right]
&\leq& C \,  \E^{Q}\Big[\,|\delta \xi|^{2p} +(\int_0^T |\delta l(s,Y_s, Z_s)|\ud s)^{2p} \Big].
\ee
Moreover we have
\ben\label{apriori2}
\E^Q[ \int_0^T |\delta Y_t|^2\ud s]+\E^{Q} \big[ \int_0^T |\delta Z_s|^2\ud s \big]
 &\leq &C \E^{Q}\left[\,|\delta \xi|^{2} +\Big(\int_0^T |\delta l(s,Y_s, Z_s)|\ud s\Big)^{2}\right].
\een
\end{lemma}

\begin{proof}[Proof of Theorem \ref{apriori-quadra}]
The arguments are similar to those of the proof of Theorem \ref{apriori.estimate}. Just make use of Lemma \ref{bmoeigen}.
\end{proof}

\section{Proof of the differentiability}\label{sec:proof:diff}
We now approach the problem of differentiability of the solutions
of a quadratic BSDEs with respect to a vector parameter on which
the terminal condition depends differentiably. We start with the
proof of the weaker property of Theorem \ref{paradiffbsde2}. Our
line of reasoning will be somewhat different from the one used for
instance by Kunita \cite{kunita} in the proof of the
diffeomorphism property of smooth flows of solutions of stochastic
differential equations. He starts with formally differentiating
the stochastic differential equation, and showing that the
resulting equation possesses a solution. The latter is then used
explicitly in moment estimates for \emph{its deviation from
difference quotients} of the original equation. The estimates are
then used to prove pathwise convergence of the difference
quotients to the solution of the differentiated SDE. We emphasize
that in our proofs, we will have to derive moment estimates for
\emph{differences of difference quotients} instead. They will
allow us to show the existence of a derivative process in a Cauchy
sequence type argument using the completeness of underlying vector
spaces, which of course will be the solution process of the
formally differentiated BSDE. So our procedure contains the
statement of the existence of a solution of the latter as a
by-product of the proof of the Theorem \ref{paradiffbsde2}. It is
not already available as a good candidate for the derivative
process, since, as we stated earlier, the formally differentiated
BSDE is a globally Lipschitz one with random Lipschitz constants
for which the classical existence theorems do not immediately
apply. Throughout assume that $f(t,x,y, z) = l(t,x,y,z) + \alpha
|z|^2$ and $\xi(x)$ satisfy (C\ref{para1}) and (C\ref{para2})
respectively.

For all $x \in \R^n$ let $(Y^x_t, Z^x_t)$ be a solution of the BSDE
(\ref{parabsde}). It is known that the solution is unique and that
$(Y^x, Z^x) \in \cR^\infty(\R^1) \times \Lzwei(\R^d)$ (see
\cite{00Kob}).

It follows from Lemma 1 in \cite{06Mor} that there exists a
constant $D > 0$ such that for all $x \in \R^n$ we have $\|(Z^x
\cdot W)_T \|_{BMO_2} \le D$. Now let $r >1$ be such that $\Psi(r) >
2 \alpha D$ (see property 3) of Lemma \ref{bmoeigen}), and denote as
before by $q$ the conjugate exponent of $r$.

\begin{proof}[Proof of Theorem \ref{paradiffbsde2}]
To simplify notation we assume that $M > 0$ is a constant such
that $\xi(x)$, $x \in \R^n$, and the derivatives of $l$ in
$(y,z)$ are all bounded by $M$. We first show that all the
partial derivatives of $Y$ and $Z$ exist. Let $x\in \R^n$ and $e_i
= (0, \ldots, 1, \ldots 0)$ be the unit vector in $\R^n$ the $i$th
component of which is $1$ and all the others $0$. For all $h \not=
0$, let $U^{h}_t = \frac1h (Y^{x+e_i h}_t - Y^x_t)$, $V^{h}_t =
\frac1h (Z^{x+he_i}_t - Z^x_t)$ and $\zeta^{h} = \frac1h (\xi(x+h
e_i) - \xi(x))$.

Let $p > 1$. Note that for all $h\not=0$ \be U^h_t = \zeta^h -
\int_t^T V^h_s \ud W_s + \int_t^T \frac1h [f(s,x+h e_i,Y^{x+h
e_i}_s,Z^{x+h e_i}_s)-f(s,x,Y^x_s,Z^x_s)] \ud s. \ee To simplify
the last term we use a line integral transformation. For all
$(\om, t) \in \Omega  \times \R_+$ let $s_{x,h} = s_{x,h}(\om, t):
[0,1] \to \R^{n+1+d}$ be defined by $s_{x,h}(\theta)=(x+\theta h
e_i, Y^x_t+\theta(Y^{x+h e_i}_t-Y^{x}_t), Z^x_t+\theta(Z^{x+h
e_i}_t-Z^x_t))$. Though $s_{x,h}$ depends on $i$ we omit to indicate this dependence for ease of notation. Note that $\frac1h s_{x,h}'(\theta)=(e_i, U^h_t,
V^h_t)$. Moreover, $A^{x,h}_t = \int_0^1 \frac{\partial
l}{\partial x_i} (s_{x,h}(\theta)) \ud \theta$, $G^{x,h}_t =
\int_0^1 \frac{\partial l}{\partial y} (s_{x,h}(\theta)) \ud
\theta$ and $I^{x,h}_t = \int_0^1 \frac{\partial l}{\partial z}
(s_{x,h}(\theta))\ud \theta$ are $(\cF_t)$-adapted processes
satisfying
\begin{eqnarray*}
\frac1h [l(t,x+he_i,Y^{x+h e_i}_t,Z^{x+h
e_i}_t)-l(t,x,Y^x_t,Z^x_t)] &=& \int_0^1 \langle \nabla
l(s_{x,h}(\theta)), s_{x,h}'(\theta) \rangle \ud \theta\\ &=&
A^{x,h}_t+G^{x,h}_t U^h_t  + I^{x,h}_t V^h_t.
\end{eqnarray*}
Since the derivatives of $l$ are bounded by $M$, $G^{x,h}_t$ and
$I^{x,h}_t$ are bounded by $M$ as well. However, we stress that
$A_t^{x,h}$ is not necessarily bounded. We define two random
functions $m^{x,h}_s(u,v)$ and $m_s(u,v)$ from $\R^{1+d}$ to $\R$ such
that $(u,v) \mapsto m^{x,h}_s(u,v) = (G^{x,h}_s u  + I^{x,h}_s v)$
and $m_s(u,v) = \left[ \partial_y l(s,Y_s^{x},Z_s^{x}) u +
\partial_z l(s,Y_s^{x},Z_s^{x}) v \right]$. Observe that these
functions satisfy (A3) and that they are Lipschitz continuous and
differentiable in $(u,v)$. In these terms, \be U^h_t = \zeta^h -
\int_t^T V^h_s \ud W_s + \int_t^T [m^{x,h}_s(U^h_s, V^h_s) +
A_s^{x,h} + \alpha (Z^{x+he_i}_s + Z^x_s) V^h_s] \ud s, \ee and thus
we obtain an equation as modelled by (\ref{bsde-lin-001}). Notice
that for all $h$, $h'\not=0$ the pairs $(U^h,V^h)$ and
$(U^h-U^{h'}, V^h-V^{h'})$ satisfy assumptions (A4) and (A5).
Therefore Theorem \ref{comparison.theo} implies with $\delta A_t=
A_t^{x,h}-A_t^{x,h'}$ \ben \label{261006-11para}
&&\E\Big[\sup_{t\in [0,T]} |U^h_t - U^{h'}_t|^{2p}\Big]\nonumber\\
&&\le C \Big\{ \E\Big[|\zeta^h - \zeta^{h'}|^{2p q^2}+\left( \int_0^T |m^{x,h}_s(U^h_s,V^h_s) - m^{x',h'}_s(U^h_s,V^h_s)+\delta A_s| ds \right)^{2pq^2}\Big]^\frac{1}{q^2}\\
&&\quad+ \E\Big[|\zeta^{h'}|^{2pq^2}+\Big(\int_0^T |A_s^{x,h'}|\ud
s\Big)^{2pq^2}\Big]^{\frac{1}{2q^2}} \E\Big[\left(\int_0^T \alpha^2
|Z^{x+he_i}_s - Z^{x+h'e_i}_s|^2 ds \right)^{2p q^2}
\Big]^{\frac{1}{2 q^2}} \Big\}. \nonumber \een Condition
(C\ref{para2}) implies that $\E[|\zeta^h - \zeta^{h'}|^{2p q^2}]$
converges to zero as $h$, $h' \to 0$. Moreover,  for some
open set $\cO$ containing $0$ we have $\sup_{h' \in
\cO\setminus \{0\}} (\E |\zeta^{h'}|^{2pq^2}) < \infty$. Due to Condition (C\ref{para1}), we
may also assume that $\sup_{h' \in \cO\setminus \{0\}}
\E[(\int_0^T |A_s^{x,h'}| \ud s)^{2pq^2}] < \infty$. Moreover,
\[\lim_{h \to 0} \E(\int_0^T |l(s,x+he_i,Y^x_s, Z^x_s) -
l(s,x,Y^x_s, Z^x_s)|\ud s)^{\beta} = 0\] for all $\beta \ge 1$, and
therefore, with Theorem \ref{apriori-quadra}, the third summand on
the right hand side of (\ref{261006-11para}) converges to zero as
$h$, $h' \to 0$.

In order to prove convergence of the second summand let $P \otimes
\lambda$ be the product measure of $P$ and the Lebesgue measure
$\lambda$ on $[0, T]$. It follows from Theorem
\ref{apriori-quadra} that $Z^{x+he_i}$ converges to $Z^x$ in
measure relative to $P \otimes \lambda$. Moreover, for all
$t\in[0,T]$, $Y^{x+he_i}_t$ converges to $Y^x_t$ in probability.
Since the partial derivatives $l_y$ and $l_z$ are continuous and
bounded, dominated convergence implies $\lim_{h\to 0} \E^{P}\left( \int_0^T |m^{x,h}_s(U^h_s,V^h_s) - m_s(U^h_s,V^h_s)| ds
\right)^{2pq^2} = 0$. Condition (C\ref{para1}) guarantees
$\lim_{h\to 0} \E^{P} \left( \int_0^T |\delta A_s| ds
\right)^{2pq^2} = 0$, and hence, $\lim_{h, h' \to 0}
\E\Big[\sup_{t\in [0,T]} |U^h_t - U^{h'}_t|^{2p}\Big] = 0$.

Finally, Theorem \ref{comparison.theo} and an estimation similar to (\ref{261006-11para}) yield
\[\lim_{h, h'\to
0} \E\left(\int_0^T|V^h_s - V^{h'}_s|^2 \ud s\right)^p = 0.\]

Now let $(h_n)$ be a sequence in $\R\setminus\{0\}$ converging to
zero. Then, since $\cR^{2p}(\R^1)$ and $\mathbb{L}^{2p}(\R^d)$ are
Banach spaces, the sequence $U^{h_n}$ converges to a process
$\frac{\partial}{\partial x_i} Y^x_t$, and $V^{h_n}$ to a process
$\frac{\partial}{\partial x_i} Z^x_t$ with respect to the
corresponding norms. By convergence term by term for the
difference quotient version of the quadratic BSDE and its formal
derivative, which follows from our a priori estimates, we see that
the pair $(\frac{\partial}{\partial x_i}
Y^x_t,\frac{\partial}{\partial x_i} Z^x_t)$ is a solution of the
BSDE \be \frac{\partial Y_t^{x}}{\partial x_i}&=&
\frac{\partial}{\partial x_i} \xi(x) -\int_t^T \frac{\partial
Z_s^{x}}{\partial x_i}\ud W_s
 + \int_t^T \Big[ \partial_{x_i} l(s,x,Y_s^{x},Z_s^{x})  \\
&&\qquad + \partial_y l(s,x,Y_s^{x},Z_s^{x}) \frac{\partial
Y_s^{x}}{\partial x_i} + \partial_z
l(s,x,Y_s^{x},Z_s^{x})\frac{\partial Z_s^{x}}{\partial x_i}  + 2
\alpha Z^x_s \frac{\partial Z_s^{x}}{\partial x_i}\Big] \ud s. \ee
Similarly to the first part one can show that $\lim_{h\to
0}\E\Big[\sup_{t\in [0,T]} |U^h_t - \frac{\partial}{\partial x_i}
Y^x_t|^{2p}\Big] = 0$ and $\lim_{h\to 0} \E\left(\int_0^T|V^h_s -
\frac{\partial}{\partial x_i}Z^x_s|^2 \ud s\right)^p = 0$, and
thus $\R^n \to \cR^{2p}(\R^1) \times \mathbb{L}^{2p}(\R^d)$, $x
\mapsto (Y^x_t, Z^x_t)$ is partially differentiable. The a priori
estimates of Theorem \ref{comparison.theo} imply that the mapping
$x \mapsto (\nabla Y^x_t, \nabla Z^x_t)$ is continuous and hence,
$(Y^x_t, Z^x_t)$ is totally differentiable. Since
differentiability with respect to $2p$th moments implies
differentiability with respect to all inferior moments above 1, we
have established the result.
\end{proof}
As a byproduct of the previous proof we obtain that for every
$x\in \R^n$ there exists a solution $(\nabla Y^x_t,\nabla Z^x_t)$
of the BSDE (\ref{pararealdiff}).\par\medskip

We now proceed with the proof of Theorem \ref{paradiffbsde}, in
which we claim pathwise continuous differentiability. To be
consistent with the previous proof, we will again compare
difference quotients varying in $h$. To this end we need the
following estimates.
\begin{lemma}\label{parahesti}
Suppose (C3) is satisfied and that $l$ and the derivatives of $l$
are all Lipschitz continuous in $(x,y,z)$. Then for all $p
>1$ there exists a constant $C> 0$, dependent only on $p$, $T$, $M$
and $D$, such that for all $x, x'\in\R^n$, $h, h' \in \R$ and $i
\in\{1, \ldots, n\}$, 
\be
& & \E \left[ \sup_{t\in[0,T]}\Big|Y_t^{x+he_i}-Y_t^{x'+h'e_i}\Big|^{2p} + \left( \int_0^T |Z^{x+he_i}_s - Z^{x'+h'e_i}_s|^2 ds
\right)^{p} \right] \\ && \qquad \qquad\qquad \qquad\le C \big(|x-x'|^2+|h-h'|^2\big)^p.
\ee
\end{lemma}
\begin{proof}
This follows from Theorem \ref{apriori-quadra}, where we put
$l_1(s,y,z) = l(s, x+h e_i, y,z),$ $l_2(s,y,z) = l(s, x'+h' e_i,
y,z).$
\end{proof}
The preceding Lemma immediately implies a first pathwise
smoothness result in $x$ for the process $Y^x$. In fact,
Kolmogorov's continuity criterion applies and yields a
modification of $Y^x$ which is continuous in $x$. More precisely:
\begin{corollary}\label{paracontiy}
There exists a process $\hat Y^x$ such that for all $(t,
\om)\in[0,T]\times \Omega$, the function $x \mapsto \hat
Y^x_t(\om)$ is continuous, and for all $(t,x)$ we have $\hat Y^x_t
=Y^x_t$ almost surely.
\end{corollary}
Let $e_i$ be a unit vector in $\R^n$. For all $x\in \R^n$ and $h
\not= 0$, let $U^{x,h}_t = \frac1h (Y^{x+h e_i}_t - Y^x_t)$,
$V^{x,h}_t = \frac1h (Z^{x+h e_i}_t - Z^x_t)$ and $\zeta^{x,h} =
\frac1h \big(\xi(x+h e_i) - \xi(x)\big)$. If $h = 0$, then define
$U^{x,0}_t = \frac{\partial}{\partial x_i} Y^x$, $V^{x,0}_t =
\frac{\partial}{\partial x_i} Z^x$ and $\zeta^{x,0} =
\frac{\partial}{\partial x_i} \xi(x)$. The proof of Theorem
\ref{paradiffbsde} will be based on the following result on the
usual difference of difference quotients. Knowing a ''good
candidate`` for the derivative from Theorem \ref{paradiffbsde2} we
allow $h=0$ this time, by replacing the difference quotient with
this candidate.
\begin{lemma}\label{201206-2}
Let $p > 1$ and $\cO \subset \R^{n+1}$ be an open set contained in
a ball of radius $\kappa$. Suppose that Condition (C3) holds and
that $l$ and the derivatives of $l$ in $(x,y,z)$ are Lipschitz
continuous in $(x,y,z)$ with Lipschitz constant $L>0$. Then there
exists a constant $C$, depending on $\kappa$, $L$, $p$, $T$, $M$,
$D$, such that for all $(x,h)$ and $(x',h') \in \cO$, \ben
\label{parakolmo1} \E \Big[\sup_{t\in [0,T]} | U^{x,h}_t -
U^{x',h'}_t|^{2p}\Big] \le C (|x-x'|^2+|h - h'|^2)^{p}. \een
\end{lemma}
\begin{proof}
Throughout the proof, $C_1, C_2, \ldots$ are constants depending
on $\kappa$, $L$, $p$, $T$, $M$, $D$.

Since $\cO$ is bounded, (C3) implies that for every $r > 1$ there
exists a constant $C_1$ such that for all $(x,h) \in \cO$ we have
$\E(\sup_{t\in [0,T]} |\zeta^{x,h}_t|^{2r}) < C_1$. Now let $s_{x,h}$,
$m^{x,h}$, $A^{x,h}$, $G^{x,h}$, $I^{x,h}$ and $U^{x,h}$ be
defined as in the proof of Theorem \ref{paradiffbsde}, and denote
$A^{x,0} = \frac{\partial l}{\partial x}(x,Y^x, Z^x)$, $G^{x,0} =
\frac{\partial l}{\partial y}(x,Y^x, Z^x)$, etc. Then the estimate
(\ref{parakolmo1}) will be deduced from the inequality
\ben\label{201206-1}
&&\E\Big[\sup_{t\in [0,T]} |U^{x,h}_t - U^{x',h'}_t|^{2p}\Big]\\
&\le& C_2 \Big\{ \E\Big[|\zeta^{x,h} - \zeta^{x',h'}|^{2p q^2}\Big]^\frac{1}{q^2}\nonumber\\
&+& \E\Big[ \Big( \int_0^T |m^{x',h'}_s(U_s^{x,h},V_s^{x,h}) - m^{x,h}_s(U_s^{x,h},V_s^{x,h})|+|A_s^{x,h}-A_s^{x',h'}| ds \Big)^{2pq^2}\Big]^\frac{1}{q^2} \nonumber\\
&+&\E\Big[|\zeta^{x',h'}|^{2pq^2}+\big(\int_0^T |A^{x',h'}_s|\ud
s\big)^{2pq^2}\Big]^{\frac{1}{2q^2}} \E\Big[  \left( \int_0^T
\alpha^2 |Z^{x'+h'e_i}_s - Z^{x+he_i}_s|^2 ds \right)^{2p q^2}
\Big]^{\frac{1}{2 q^2}}\Big\} \nonumber \een which follows from
Theorem \ref{comparison.theo}. We first analyze the order of the
convergence of
\[B_1(x,x',h,h')= \E \Big[\left( \int_0^T |m^{x',h'}_s(U_s^{x,h},V_s^{x,h}) - m^{x,h}_s(U_s^{x,h},V_s^{x,h})| \ud s \right)^{2pq^2} \Big]^\frac{1}{q^2}
\textrm{ as }h, h' \to 0.\] To this end notice that \be
B_1(x,x',h,h')
&\le& C_3 \Big\{ \left(\E \left( \int_0^T |G^{x',h'}_t - G^{x,h}_t| |U^{x,h}_t| \ud t \right)^{2pq^2}\right)^\frac{1}{q^2}\\
&&\qquad\qquad + \left(\E \left( \int_0^T |I^{x',h'}_t - I^{x,h}_t| |V^{x,h}_t| \ud t \right)^{2pq^2}\right)^\frac{1}{q^2}\Big\}.
\ee
Then
\be
&&\int_0^T |G^{x',h'}_t - G^{x,h}_t| |U^{x,h}_t| \ud t\\
&&\quad\le \sup_{t\in[0,T]} |U^{x,h}_t| \int_0^T |G^{x',h'}_t - G^{x,h}_t| \ud t \\
&&\quad\le \sup_{t\in[0,T]} |U^{x,h}_t| \int_0^T \left( \int_0^1 |\partial_y l (s_{x',h'}(\theta)) - \partial_y l(s_{x,h}(\theta))| \ud \theta \right) \ud t \\
&&\quad\le  \sup_{t\in[0,T]} |U^{x,h}_t| \int_0^T \int_0^1 L|s_{x',h'}(\theta) - s_{x,h}(\theta)| \ud \theta \ud t \\
&&\quad\le C_4 \sup_{t\in[0,T]} |U^{x,h}_t| \Big( |x'-x|+|h'-h|+\sup_{t\in[0,T]} |Y^{x'}_t - Y^x_t| + \sup_{t\in[0,T]}|Y^{x'+h' e_i}_t - Y^{x+he_i}_t| \\
&&\quad\qquad+ \int_0^T  (|Z^{x'}_t - Z^x_t| + |Z^{x'+h' e_i}_t - Z^{x+he_i}_t|) \ud t \Big),
\ee
and, by applying H\"older's inequality we obtain with Lemma \ref{parahesti}
\be
\left( \E \left( \int_0^T |G^{x',h'}_t - G^{x,h}_t| |U^{x,h}_t| \ud t \right)^{2pq^2} \right)^\frac{1}{q^2} &\le& C_5 \left(|h-h'|^{2} + |x -x'|^{2}\right)^{p}.
\ee
Similarly, $\left(\E \left[ \int_0^T |I^{x',h'}_t - I^{x,h}_t| |V^{x,h}_t| \ud t \right]^{2pq^2} \right)^\frac{1}{q^2} \le C_6 \left(|h-h'|^{2} + |x -x'|^{2}\right)^{p}$, and so we conclude $B_1(x,x',h,h') \le C_7 \left(|h-h'|^{2} + |x -x'|^{2}\right)^{p}$.

By using similar arguments we get \be \E \left[\Big(
\int_0^T |A_t^{x,h}-A_t^{x',h'}|\ud t
\Big)^{2pq^2}\right]^\frac{1}{q^2} &\le&
\E\left[ \left(\int_0^T \int_0^1 |\partial_x l (s_{x',h'}(\theta)) - \partial_x l(s_{x,h}(\theta))| \ud \theta  \ud t \right)^{2pq^2}\right]^\frac{1}{q^2} \\
&\leq&
C_8 \E\left[ \left(\int_0^T \int_0^1 |s_{x',h'}(\theta) - s_{x,h}(\theta)| \ud \theta  \ud t \right)^{2pq^2}\right]^\frac{1}{q^2} \\
&\le& C_9 \left(|h-h'|^{2} + |x -x'|^{2}\right)^{p}.
\ee
Theorem \ref{apriori-quadra} and the Lipschitz continuity of $l$ imply
\be
&&\E\Big[ \left( \int_t^T |Z^{x+he_i}_s - Z^{x'+h' e_i}_s|^2 ds \right)^{2p q^2}\Big]^\frac{1}{2q^2}\\
&\le& C_{10} \E\Big[|(\xi(x+he_i) - \xi(x'+h' e_i)|^{4p q^4}\\
&&\qquad + \Big(\int_0^T|l(s,x+he_i,Y^{x+he_i},Z^{x+he_i})-l(s,x'+h'e_i,Y^{x+he_i},Z^{x+he_i})|\ud s\Big)^{4pq^4}\Big]^\frac{1}{2q^4}\\
&\le& C_{11} (|x-x'|^2 + |h-h'|^2)^{p}. \ee Finally, (C3) yields
$(\E|\zeta^{x,h} - \zeta^{x',h'}|^{2p q^2})^\frac{1}{q^2} \le
C_{12} (|x-x'|^2 + |h-h'|^2)^{p}$, and hence \be \E
\Big[\sup_{t\in [0,T]} | U^{x,h}_t - U^{x',h'}_t|^{2p}\Big] \le
C_{13} (|x-x'|^2 + |h-h'|^2)^{p}. \ee
\end{proof}
\begin{proof}[ Proof of Theorem \ref{paradiffbsde}]
To simplify notation we may assume that (\ref{parakolmo1}) is
satisfied for $\cO = \R^{n+1}$. Assume that $Y^x_t$ is continuous
in $x$ (see Corollary \ref{paracontiy}). Lemma \ref{201206-2} and
Kolmogorov's continuity criterion imply that  $U^{x,h}_t$ has a
modification $\hat U^{x,h}_t$ continuous in $(x,h)$. Define
$\frac{\partial}{\partial x_i} Y^x_t = \hat U^{x,0}_t$ and note that
we obtain thus a continuous version of the solution of the BSDE
(\ref{pararealdiff}). For all $(x,h) \in \Q^{n+1}$ let $N(x,h)$ be
a null set such that for all $\om \notin N(x,h)$ we have $\hat
U^{x,h}_t(\om)=U^{x,h}_t(\om)$. Then, $N = \bigcup_{(x,h) \in
\Q^{n+1}} N(x,h)$ is a null set such that for all $\om \notin N$
the following implication holds: If $q_k \in \Q^n$ and $r_k \in
\Q\setminus \{0\}$ are sequences with $\lim_{k \to \infty} q_k = x
\in \R^n$ and $\lim_{k\to \infty}r_k = 0$, then
\[ \lim_{k \to \infty} \frac{1}{r_k}(Y^{q_k+r_k e_i}_t -Y^{q_k}_t) = \frac{\partial}{\partial x_i} Y^x_t. \]
As a consequence of this and the subsequent Lemma
\ref{rational.para}, $Y^x_t(\om)$ is continuously partially
differentiable relative to $x_i$ if $\om \notin N$. Since we can
choose such a null set for any $i \in\{1, \ldots, n\}$, total
differentiability follows and the proof is complete.
\end{proof}
\begin{lemma}\label{rational.para}
Let $f:\R^n\to \R$ be a continuous function and $g:\R^n \to \R^n$
a continuous vector field. Suppose that for all sequences $q_k \in
\Q^n$ with $q_k \to x \in \R^n$ and $r_k \in \Q\setminus \{0\}$
with $r_k \to 0$ we have
\[ \lim_{k \to \infty} \frac{1}{r_k} (f(q_k + r_k e_i) - f(q_k)) = g_i(x), \]
where $1 \le i \le n$. Then $f$ is differentiable and $\nabla f = g$.
\end{lemma}
\begin{proof}
To simplify notation assume that $n = 1$. Let $x_k \in \R$ with
$x_k \to x \in \R$ and $h_k \in \R\setminus \{0\}$ with $h_k \to
0$. Since $f$ is continuous we may choose $q_k \in \Q$ and $r_k
\in \Q \setminus \{0\}$ such that $|f(q_k) - f(x_k)| \le
\frac{|h_k|}{2^k}$, $|f(q_k + r_k) - f(x_k + h_k)| \le
\frac{|h_k|}{2^k}$ and $|\frac{1}{r_k} - \frac{1}{h_k}| \le
\frac{1}{2^k}$. Then \be
&&\!\!\!\!\!\!\!\!\!\!\!\! \!\!\!\!\!\!\!\!\!\!\!\! |\frac{1}{h_k} (f(x_k + h_k) - f(x_k)) - g(x)|\\
\quad&\le& \Big|\frac{1}{h_k}[\big(f(x_k + h_k) - f(x_k)\big) - \big(f(q_k + r_k) - f(q_k)\big)]\Big| \\
\quad& & \quad + |(\frac{1}{h_k}-\frac{1}{r_k}) (f(q_k + r_k) - f(q_k))|
+ |\frac{1}{r_k} \big(f(q_k + r_k) - f(q_k)\big) - g(x)| \\
\quad&\le& 2 \frac{1}{2^k} + \frac{1}{2^k} |f(q_k + r_k) - f(q_k)| + |\frac{1}{r_k} (f(q_k + r_k) - f(q_k)) - g(x)| \\
\quad&\to& 0, \qquad (k \to \infty), \ee and hence $f$ is
partially differentiable. Since the partial derivatives $g_i$ are
continuous, $f$ is also totally differentiable.
\end{proof}

\section{Differentiability of quadratic Forward-Backward
SDEs}\label{sec:forward-backward}

In this section we will specify the results obtained in the
preceding sections to BSDEs where the terminal conditions
are determined by a forward SDE driven by the same Brownian motion
as the BSDE. When considering BSDEs with terminal condition
determined by a forward SDE we will need regularity of the forward
equation. This will be guaranteed if the coefficients are
functions belonging to the following space.

Throughout this section let again $n$ be a positive integer and
$W$ a $d$-dimensional Brownian motion.
\begin{defi} Let $k$, $m \ge 1$. We denote by $\bf B^{k \times m}$ the set of all functions $h:[0,T] \times \R^n \to \R^{k \times m}$,
$(t,x) \mapsto h(t,x)$, differentiable in $x$,  for which there
exists a constant $C >0$ such that $\sup_{(t,x)\in[0,T]\times\R^n}
\sum_{i=1}^n \left| \frac{\partial h(t,x)}{\partial x_i} \right|
\le C$; and for all $t\in[0,T]$ we have $\sup_{x \in\R^n}
\frac{|h(t,x)|}{1+|x|} \le C$ and $x \mapsto \frac{\partial
h(t,x)}{\partial x}$ is Lipschitz continuous with Lipschitz
constant $C$.

With any pair $h \in \bf B^{n \times 1}$ and $\sigma \in \bf
B^{n\times d}$ we associate the second order differential operator
$\cL = \sum_{i=1}^n h_i(\cdot) \frac{\partial}{\partial x_i} +
\frac12 \sum_{i, j=1}^n [\sigma\sigma^T]_{ij}(\cdot)
\frac{\partial^2}{\partial x_i \partial x_j}$.
\end{defi}

We will consider Forward-Backward SDEs (FBSDEs) of the form
\begin{equation}\label{fbsde}
\left\{
\begin{array}{lll} X_t^x&=&x+\int_0^t b(s,X_s^{x})\ud s+\int_0^t \sigma(s,X_s^{x})\ud W_s,\qquad x\in\R^n, \\
Y_t^{x}&=&g(X_T^{x})+\int_t^T f(s,X^x_s, Y_s^x,Z_s^{x})\ud
s-\int_t^T Z_s^{x}\ud W_s,
\end{array}
\right.
\end{equation}
where the coefficients satisfy the following assumptions:
\begin{enumerate}[({D}1)]
\item\label{fbforward}
$\sigma\in \bf B^{n \times d}$, $b \in \bf B^{n\times 1}$,
\item\label{fbbackward}
$f:\Omega\times[0,T]\times\R^n\times\R\times\R^d\to\R$ and $g:\R^n\to \R$ are measurable functions such that
$f(\omega,t,x,y,z)=l(\omega,t,x,y,z)+\alpha |z|^2$, where $l(\omega, t, x, y, z)$ is globally Lipschitz and continuously differentiable in $(x,y,z)$,
\item\label{fbterminal}
$g:\R^n\to \R$ is a twice differentiable function such that $\nabla g \cdot \sigma \in \bf B^{1 \times d}$ and $\cL g \in \bf B^{1\times 1}$.
\end{enumerate}
It follows from standard results on SDEs and from Theorem 2.3 in
\cite{00Kob} that there exists a solution $(X^x,Y^x,Z^x)$ of
Equation (\ref{fbsde}). As we will show, the results of Section
\ref{diffsection} imply that $(X^x,Y^x,Z^x)$ is differentiable in
$x$ and that the derivatives $(\nabla X^x,\nabla Y^x,\nabla Z^x)$
solve the FBSDE
\begin{equation}
\begin{array}{lll}
\nabla X_t^{x}&=&1+\int_0^t \partial_x b(X_s^{x})\nabla X^x_s\ud s
  +\int_0^t \partial_x\sigma(X_s^{x})\nabla X^x_s\ud W_s,\\
\nabla Y_t^{x}&=& \partial_x g(X_T^{x})\nabla X^x_T-\int_t^T \nabla Z_s^{x}\ud W_s+\int_t^T \Big[ \partial_x l(s,X^x_s,Y_s^{x},Z_s^{x}) \nabla
X_s^{x} + \partial_y l(s,X^x_s,Y_s^{x},Z_s^{x}) \nabla
Y_s^{x} \\
&&  + \partial_z l(s,X^x_s,Y_s^{x},Z_s^{x}) \nabla Z_s^{x} + 2
\alpha Z^x_s \nabla Z_s^{x}\Big] \ud s.\label{realdiff}
\end{array}
\end{equation}
Our first result parallels Theorem \ref{paradiffbsde2} in which
differentiability with respect to vector space topologies is
treated.
\begin{theo}\label{diffbsde2}
Let (D\ref{fbforward}) and (D\ref{fbbackward}) be satisfied and
assume that $g: \R^n \to \R$ is bounded and differentiable.
Moreover, suppose that $\frac{\partial l}{\partial x}(t,x,y,z)$ is
Lipschitz continuous in $x$. Then for all $p \ge 2$, the function
$\R^n \to \cR^p(\R^n) \times \cR^p(\R^1) \times \Lp(\R^d)$, $x
\mapsto (X^x, Y^x, Z^{x})$, is differentiable, and the derivative
is a solution of the BSDE (\ref{realdiff}).
\end{theo}
\begin{proof}
By standard results, the mapping $\R^n \to \cR^p(\R^1)$, $x
\mapsto X^x$ has a continuous version (which we assume being
identical to the given one), and for all $p
> 1$ there exists a constant $C \in \R_+$ such that for $x$, $x'
\in \R^n$ we have $\E(|X^x -X^{x'}_t|^{2p}) \le C |x - x'|^{2p}.$
See for example Lemma 4.5.4 and Lemma 4.5.6 in \cite{kunita}. In
order to be able to apply Theorem \ref{paradiffbsde2}, we need to
verify Condition (C\ref{para1}). For this purpose, note that \be
&&\E\left[\left( \int_0^T|l(t,X^{x}_t, Y^x_t,Z^x_t)-l(t,X^{x'}_t, Y^x_t,Z^x_t)| \ud t\right)^{2p}\right]\\
& &\qquad\le \E\left[\left(\int_0^T L|X^{x}_t-X^{x'}_t| \ud t\right)^{2p}\right] \le \tilde{C}(L,T,C) |x-x'|^{2p},
\ee
where $L$ is a Lipschitz constant of $l$. Similarly,
\[
\E\big[\left(\int_0^T |\frac{\partial}{\partial x}l(t,X^{x}, Y^x,Z^x)-\frac{\partial}{\partial x}l(t,X^{x'}, Y^x,Z^x) \ud t|\right)^{2p}\big]
 \le C' |x-x'|^{2p}.
\] This proves (C\ref{para1}). Moreover, notice that $\xi(x) = g(X^x_T)$ satisfies Condition (C\ref{para2}).
Thus the statement follows from Theorem \ref{paradiffbsde2}.
\end{proof}
If in addition Condition (D\ref{fbterminal}) is satisfied, we
again obtain a sharper result stating pathwise continuous
differentiability of an appropriate modification of the solution
process.
\begin{theo}\label{fbdiffbsde}
Assume that $(D\ref{fbforward})$, $(D\ref{fbbackward})$ and
$(D\ref{fbterminal})$, and suppose that the partial derivatives of
$l$ in the variables $(x,y,z)$ are Lipschitz continuous. Then
there exists a function $\Omega \times [0, T] \times \R^n \to
\R^{n+1+d}$, $(\om, t, x) \mapsto (X^x_t, Y^{x}_t, Z^x_t)(\om)$,
such that for almost all $\om$, $X^x_t$ and $Y^x_t$ are continuous
in $t$ and continuously differentiable in $x$, and for all $x$,
$(X^x_t, Y^x_t, Z^x_t)$ is a solution of (\ref{fbsde}).
\end{theo}
Let $M > 0$ be a constant such that $g$, the derivatives of $g$,
$b$ and $\sigma$, and the partial derivatives of $l$ in $(x,y,z)$
are all bounded by $M$. For all $x \in \R$ let $(X^x_t, Y^x_t,
Z^x_t)$ be the solution of the FBSDE (\ref{fbsde}). To correspond
formally to Theorem \ref{paradiffbsde}, in the setting of our
FBSDE we have to work with $$\tilde{l}(\omega, t, x, y, z) =
l(\omega, t, X^x_t(\omega), y,z).$$ But this functional fails to
be globally Lipschitz in $x$. This is why we have to modify
slightly the proof of Theorem \ref{paradiffbsde}, and cannot just
quote it. We start by showing that $\xi(x) = g(X^x_T)$ satisfies
Condition (C3).
\begin{lemma}\label{c3bed}
For all $x\in \R^n$, $h \not= 0$ and $i \in \{1, \ldots, n\}$, let
$\zeta^{x,h,i} = \frac1h (g(X^{x+h e_i}_T) - g(X^x_T))$. Then for
every $p > 1$ there exists a $C>0$, dependent only on $p$ and $M$,
such that for all $x$, $x' \in \R^n$ and $h$, $h' \not= 0$, \be
\E\Big[|\zeta^{x,h,i} - \zeta^{x',h',i}|^{2p} \Big] \le C
(|x-x'|^2 +|h - h'|^2)^p. \ee Moreover, for all $t \in[0,T]$,
\be \E\Big[|\zeta^{x,h,i} - g'(X^{x}_T) \frac{\partial}{\partial
x_i}X^x_T|^{2p} \Big] \le C (|x-x'|^2)^p. \ee
\end{lemma}
\begin{proof}
Note that by Ito's formula $g(X^x_t)= g(X^x_0) + \int_0^t \nabla g(X^x_s) \cdot \sigma(s, X^x_s) dW_s + \int_0^t \cL g ds$. Thus $g(X^x_t)$ is a diffusion with coefficients $\tilde \sigma(s,x)=\nabla g(x) \cdot \sigma(s, x)$ and $\tilde b(s,x)=\sum_{i=1}^n b_i(s,x) \frac{\partial g(x)}{\partial x_i} + \frac12 \sum_{i, j=1}^n \sigma_{ij}(s,x) \frac{\partial^2 g(x)}{\partial x_i \partial x_j}$. By (D\ref{fbterminal}) we have $\tilde \sigma \in \bf B^{1\times d}$ and $\tilde b \in \bf B^{1\times 1}$. Therefore, by using standard results on stochastic flows (see Lemma 4.6.3 in \cite{kunita}), we obtain the result.
\end{proof}
\begin{proof}[ Proof of Theorem \ref{fbdiffbsde}]
First note that it is well-known that $X^x$ may be chosen to be
continuous in $t$ and continuously differentiable in $x$ (see for
example Theorem 39, Ch. V, \cite{protter}). In order to prove that
$Y^x$ has such a modification as well, note that Lemma \ref{c3bed}
implies that $\xi(x) = g(X^x_T)$ satisfies Condition (C3). Now let
again $U^{x,h}_t = \frac1h (Y^{x+h e_i}_t - Y^x_t)$, $V^{x,h}_t =
\frac1h (Z^{x+h e_i}_t - Z^x_t)$ for all $x\in \R^n$ and $h \not=
0$. If $h = 0$, then define $U^{x,0}_t = \frac{\partial}{\partial
x_i} Y^x$, $V^{x,0}_t = \frac{\partial}{\partial x_i} Z^x$ and
$\zeta^{x,0} = \frac{\partial}{\partial x_i} g(X^{x}_T)$. It is
enough to show that for all open bounded sets $\cO \subset
\R^{n+1}$ there exists a constant $C$ such that for all $(x,h) \in
\cO$
\[
\E \Big[\sup_{t\in [0,T]} | U^{x,h}_t - U^{x',h'}_t|^{2p}\Big] \le C (|x-x'|^2+|h - h'|^2)^{p}.
\]
As in Lemma \ref{201206-2} we will derive this estimate from
Inequality (\ref{201206-1}). Notice that the assumptions of
Theorem \ref{fbdiffbsde} guarantee that all the terms appearing in
(\ref{201206-1}), satisfy the same properties and thus provide the
same estimates. There is one essential difference which is due to
the appearance of $X^x$ instead of $x$ in the first component of
the line described by the integral
$\gamma_{x,h}(\theta)=(X^x+\theta(X^{x+h e_i} - X^x),
Y^x_t+\theta(Y^{x+h e_i}_t-Y^{x}_t), Z^x_t+\theta(Z^{x+h
e_i}_t-Z^x_t))$. We therefore have to replace the prior
$A^{x,h}_t$ by $A^{x,h}_t = B^{x,h}_t \frac1h (X^{x+he_i}- X^x)$
with $B^{x,h}_t = \int_0^1 \frac{\partial l}{\partial x}
(\gamma_{x,h}(\theta)) \ud \theta$. Therefore we only need to show
for $(x,h) \in \cO$ \ben\label{211206-1}
\E\big[\big(\int_0^T|A_s^{x,h}-A_s^{x',h'}| \ud s
\big)^{2pq^2}\big]^\frac{1}{q^2} \le c(|x-x'|^2+|h-h'|^2)^p. \een
In fact, with $\Delta^{x,h}=\frac1h (X^{x+he_i}- X^x)$ if
$h\not=0$ and $\Delta^{x,0}=\nabla X^x$, we have \be
&&\E\big(\int_0^T|A_s^{x,h}-A_s^{x',h'}| \ud s \big)^{2pq^2} \\
&\le&\int_0^T  |\Delta^{x,h}_s| |B_s^{x,h}-B_s^{x',h'}| \ud s + \int_0^T |\Delta^{x,h}_s-\Delta^{x',h'}_s||B_s^{x',h'}| \ud s.
\ee
The first summand satisfies
\be
&&\E\big(\int_0^T |\Delta^{x,h}_s||B_s^{x,h}-B_s^{x',h'}| \ud s \big)^{2pq^2}\\
&&\qquad\le \E\big(\int_0^T |\Delta^{x,h}_s| \int_0^1  |\frac{\partial l}{\partial x} (\gamma_{x,h}(\theta)) -\frac{\partial l}{\partial x} (\gamma_{x',h'}(\theta))| \ud \theta  \ud s \big)^{2pq^2} \\
&&\qquad\le \left( \E\big(\int_0^T |\Delta^{x,h}_s| \ud s
\big)^{4pq^2}\right)^\frac12 \left(\E\big(\int_0^T \int_0^1
|\gamma_{x,h}(\theta)) -\gamma_{x',h'}(\theta)| \ud \theta \ud
s\big)^{4pq^2}\right)^\frac12. \ee Lemma 4.6.3 in \cite{kunita}
implies $\sup_{(x,h) \in \cO} \E\big(\int_0^T |\Delta^{x,h}_s| \ud
s \big)^{4pq^2} < \infty$.
Besides, \be
&& \E\big(\int_0^T \int_0^1  |\gamma_{x,h}(\theta)) -\gamma_{x',h'}(\theta)| \ud \theta \ud s\big)^{4pq^2} \\
&&\qquad\le C \, \E \Big( \sup_{t\in[0,T]} |X^{x'}_t-X^{x}_t|+ \sup_{t\in[0,T]} |X^{x'+h' e_i}_t - X^{x + h e_i}_t|
+\sup_{t\in[0,T]} |Y^{x'}_t - Y^x_t|\\&&
\qquad+ \sup_{t\in[0,T]}|Y^{x'+h' e_i}_t - Y^{x+he_i}_t|
 + \int_0^T  (|Z^{x'}_t - Z^x_t| + |Z^{x'+h' e_i}_t -
Z^{x+he_i}_t|) \ud t \Big)^{4pq^2}. \ee From this we can easily
deduce $\E\left[\big(\int_0^T
|\Delta^{x,h}_s||B_s^{x,h}-B_s^{x',h'}| \ud s
\big)^{2pq^2}\right]^\frac{1}{q^2} \le C (|x-x'|^2+|h -
h'|^2)^{p}$. Similarly, $\E\left[\big(\int_0^T
|\Delta^{x,h}_s-\Delta^{x',h'}_s||B_s^{x',h'}| \ud s
\big)^{2pq^2}\right]^\frac{1}{q^2} \le C (|x-x'|^2+|h -
h'|^2)^{p}$, hence (\ref{211206-1}) follows and the proof is
finished.
\end{proof}

\section{Malliavin differentiability of quadratic
BSDEs}\label{sec:malliavin}

In this section we shall ask for a different type of smoothness for
solutions of quadratic BSDEs, namely differentiability in the
variational sense or in the sense of Malliavin's calculus. Of
course, this will imply smoothness of the terminal condition in the
same sense. If the terminal condition is given by a smooth function
of the terminal value of a forward equation, it will also involve
variational smoothness of the forward equation.\par\medskip

Let us first review some basic facts about Malliavin calculus. We
refer the reader to \cite{nualart} for a thorough treatment of the
theory and to \cite{97KPQ} for results related to BSDEs. To begin
with, let $C_b^\infty(\R^{n\times d})$ denote the set of functions
with partial derivatives of all orders defined on $\R^{n\times d}$
whose partial derivatives are bounded.

Let $\cS$ denote the space of random variables $\xi$ of the form
\[
\xi = F\Big((\int_0^T h^{1,i}_s d W^1_s)_{1\le i\le
n},\cdots,(\int_0^T h^{d,i}_s d W^d_s)_{1\le i\le n})\Big),
\] where $F\in C_b^\infty(\R^{n\times d})$, $h^1,\cdots,h^n\in L^2([0,T]; \R^d)$.
To simplify the notation assume that all $h^j$ are written as row
vectors.

If $\xi\in \cS$ of the above form, we define the $d$-dimensional
operator $D = (D^1,\cdots, D^d):\cS\to L^2(\Omega\times[0,T])^d$ by
\[
D^i_\theta \xi = \sum_{j=1}^n \frac{\partial F}{\partial x_{i,j}}
\Big( \int_0^T h^1_tdW_t,\ldots,\int_0^T
h^n_tdW_t\Big)h^{i,j}_\theta,\quad 0\leq \theta\leq T,\quad 1\le
i\le d.
\] For $\xi\in\cS$ and $p>1$, we define the norm
\[
\lVert \xi \lVert_{1,p} = \Big(\E\Big[\,|\xi|^p+\Big(\int_0^T |D_\theta \xi|^2 \ud \theta \Big)^{\frac{p}{2}}\Big]\Big)^{\frac1p}.
\]
It can be shown (see for example \cite{nualart}) that the operator
$D$ has a closed extension to the space $\mathbb{D}^{1,p}$, the closure of
$\cS$ with respect to the norm $\lVert \cdot \lVert_{1,p}$. Observe
that if $\xi$ is $\F_t-$measurable then $D_\theta \xi=0$ for $\theta
\in (t,T]$.

We shall also consider $n$-dimensional processes depending on a time
variable. We define the space $\mathbb{L}^a_{1,p}(\R^n)$ to be the
set of $\R^n-$valued progressively measurable processes
$u(t,\omega)_{ t\in[0,T], \omega\in \Omega}$ such that
\begin{itemize}
\item[i)]
For a.a. $t\in[0,T]$, $u(t,\cdot)\in(\mathbb{D}^{1,p})^n$;
\item[ii)] $(t,\omega)\to D_\theta u(t,\omega)\in (L^2([0,T]))^{d\times n}$ admits a progressively measurable
version;
\item[iii)] $\lVert u \lVert^a_{1,p}=\E[\Big(\int_0^T |u(t)|^2\ud t\Big)^{\frac{p}{2}}+
\Big(\int_0^T\int_0^T|D_\theta u(t)|^2\ud \theta \ud
t\Big)^{\frac{p}{2}}]^{\frac1p}<\infty.$
\end{itemize}
Here, for $y\in\R^{d\times n}$ we use the norm $|y|^2 = \sum_{i,j}
(y_{i,j})^2.$

We also consider the space \[\mathbb{D}^{1,\infty}=\cap_{p>
1}\mathbb{D}^{1,p}.\]

We cite for completeness a result from \cite{nualart} that we will use in the next section.
\begin{lemma}[Lemma 1.2.3 in \cite{nualart}]\label{lemma.from.nualart}
Let $\{F_n,\,n\geq 1\}$ be a sequence of random variables in $\mathbb{D}^{1,2}$ that converges to $F$ in $L^2(\Omega)$ and such that
\[\sup_{n\in\N} \E[\lVert D F_n\lVert_{L^2}]<\infty.\] Then $F$ belongs to $\mathbb{D}^{1,2}$, and the sequence of derivatives $\{D F_n,\,n\geq 1\}$ converges to $D F$ in the weak topology of
$L^2(\Omega\times [0,T])$.
\end{lemma}

Let us now consider the BSDE \ben\label{linbsde.aux015} Y_t = \xi -
\int_t^T Z_s \ud W_s + \int_t^T f(s,Y_s,Z_s) ds. \een

Our assumptions on driver and terminal condition this time amount to

\begin{enumerate}[({E}1)]
\item\label{para3} $f:\Omega \times [0,T] \times\R\times\R^d\to\R$
is an adapted measurable function such that
$f(\omega,t,y,z)=l(\omega,t,y,z)+\alpha |z|^2$, where $l(\omega, t,
y, z)$ is globally Lipschitz and continuously differentiable in
$(y,z)$; for all $p>1$ we have $\E^P[(\int_0^T |l(\omega, t, 0,0)|^2\ud s)^{2p}]<\infty$;

\item\label{para5}for all $(t,y,z)$, the mapping $\Omega\to \R,$
$\omega\mapsto l(\omega, t,y,z)$ is Malliavin differentiable and
belongs to $\mathbb{L}^a_{1,p}(\R)$ for all $p>1.$

For any $(\omega,t,y,z)$ and $\theta\in [0,T]$, the (a.e. valid)
inequality holds true
\[
|D_\theta l(\omega,t,y,z))|\leq
\tilde{K}_\theta(\omega,t)+K_\theta(\omega,t)(|y|+|z|)
\] where $K_\theta$ and $\tilde{K}_\theta$
are positive adapted processes satisfying for all $p\geq 1$
\be
\E[\Big(\sup_{t\in[0,T]} \int _0^T |K_\theta
(t, \omega)|^2\ud \theta \Big)^p]<\infty\quad \textrm{and} \quad
\E[\Big(\int_0^T \int_0^T |\tilde{K}_\theta(t, \omega)|^2\ud
\theta\ud t\Big)^p]<\infty \ee

\item\label{para4} the random
variable $\xi$ is bounded and belongs to $\mathbb{D}^{1,\infty}.$
\end{enumerate}

We first consider the case where the terminal variable has no
further structural properties, such as depending on the terminal
value of a forward equation. For notational simplicity we shall
treat the case of one dimensional $z$ and Wiener process and so
may omit the superscript $i$ in $D^i$ etc. We will this time use
the typical Sobolev space approach, hidden in Lemma
\ref{lemma.from.nualart}, to describe Malliavin derivatives, which
are in fact derivatives in the distributional sense. In this
approach we shall employ an approximation of the driver of our
BSDE by a sequence of globally Lipschitz continuous ones, for
which the properties we want to derive are known.

Let us therefore introduce a family of truncated functions
starting with describing their derivatives by
\begin{displaymath}
g'_n(z)= \left\{
\begin{array}{cl}
-2n&,z<-n\\
2z&,|z|\leq n\\
2n&,z>n.
\end{array}
\right.
\end{displaymath}
Then we have $g_n(z)=z^2$ for $|z|\leq n$, $g_n(z) = 2n|z|-n^2$ for
$|z|>n$, and thus $|g_n(z)|\leq z^2$ and $g_n(z)\to z^2$ locally
uniformly on $\R$ for $n\to\infty$. A similar statement holds for
the derivative of $g_n(z)$: $|g'_n(z)|\leq 2|z|$ and $g'_n(z)\to 2z$
locally uniformly on $\R$ for $n\to\infty$.

With these truncation functions we obtain the following family of
BSDEs: \ben\label{linbsde.aux:24.jan.2007} Y^n_t = \xi - \int_t^T
Z^n_s dW_s + \int_t^T [l(s,Y^n_s,Z^n_s)+\alpha g_n(Z^n_s)] ds,\quad
n\in\N. \een From Proposition 2.4 of \cite{00Kob} we obtain that there
exists $(Y_s,Z_s)\in \cR^\infty(\R) \times \mathbb{L}^2(\R)$ such
that $Y^n_s\to Y_s$ uniformly in $[0,T]$ and $Z^n_s\to Z_s$ in
$\mathbb{L}^2(\R)$.

Since the truncated equations have Lipschitz continuous drivers,
Proposition 5.3 of \cite{97KPQ} guarantees that $(Y^n_t,Z^n_t)\in
\mathbb{D}^{1,2}\times \mathbb{D}^{1,2}$ with the following
Malliavin derivative \ben
D_\theta Y_t^n &=&0 \quad\textrm{ and }\quad D_\theta Z_t^n =0, \textrm{ if }t\in[0,\theta),\nonumber\\
D_\theta Y_t^n &=&D_\theta\xi +\I\Big[\partial_y l(Y^n_s,Z^n_s)D_\theta Y^n_s+\partial_z l(Y^n_s,Z^n_s)D_\theta Z^n_s
\nonumber\\&&+D_\theta l(s,Y^n_s,Z^n_s)
+\alpha g'_n(Z^n_s)D_\theta Z^n_s\Big]\ud s-\I D_\theta Z^n_s\ud W_s,\label{aux:1:10.jan.2007} \qquad
\textrm{ if }t\in[\theta,T]. \een

Now we aim at showing that the sequences $D Y^n$ and $D Z^n$ are
bounded in $\mathbb{D}^{1,2}$, in order to use Lemma
\ref{lemma.from.nualart}. This will be done by deriving a priori
estimates in the style of the preceding sections, this time
uniform in $n$. We therefore first show boundedness relative to
the auxiliary measures $Q_n:=\cE\big(\alpha \int g'_n(Z^n) dW
\big)\cdot P$, in the form of the following a priori inequality.

\begin{lemma}\label{lemma.aux:2:11.jan.2007}
Let $p>1.$ If the driver and terminal condition
satisfy hypotheses (E\ref{para3}), (E\ref{para5}) and
(E\ref{para4}), then the
following inequality holds for the BSDE (\ref{aux:1:10.jan.2007}):
\be &&\E^{Q_n}\Big[\Big(\sup_{t\in [0,T]} \int_0^T |D_\theta
Y^n_t|^2\ud \theta \Big)^p\Big]+ \E^{Q_n}\Big[\Big(\int_0^T\int_0^T
|D_\theta Z^n_s|^2\ud \theta\ud s\Big)^p\Big]
\\
&&\qquad\qquad\leq
C\E^{Q_n}\Big[ \Big(\int_0^T|D_\theta \xi|^2\ud \theta\Big)^p
 +\Big(\int_0^T\int_0^T |D_\theta l(s,Y^n_s,Z^n_s)|^2\ud \theta\ud s\Big)^p\Big]
\ee
\end{lemma}
\begin{proof}
We will derive these estimates by proceeding in the same fashion as
for Lemma \ref{apriori.estimate.lemma}. Again, $C_1, C_2,\ldots$ are constants depending on the
coefficients and $p$.

Applying It\^o's formula to $\ebt |D_\theta Y^n_t|^2$, using
Equation (\ref{aux:1:10.jan.2007}) and simplifying as we did in the
former sections we obtain (choosing $\beta=M^2+2M$) \ben
&&e^{\beta t} |D_\theta Y^n_t|^2 +\I e^{\beta s} (M|D_\theta Y^n_s|-|D_\theta Z^n_s|)^2\ud s\nonumber\\
&&\quad\leq e^{\beta T} |D_\theta \xi|^2-2\I e^{\beta s} D_\theta
Y^n_s D_\theta Z^n_s \ud \hat{W}_s + 2\I e^{\beta s} |D_\theta Y^n_s
D_\theta l(s,Y^n_s,Z^n_s)| \ud s,\label{aux:2:10.jan.2007} \een
where $ \hat{W}_t= W_t-\int_0^t \alpha g_n'(Z_s^n)\ud s,$
$t\in[0,T],$ is a $Q_n-$ Brownian motion.

We remark that since
$\E\sup_{t\in[0,T]} |D_\theta Y^n_t|^2<\infty$ we have
$\E\Big[(\int_0^T |D_\theta Y^n_s|^2|D_\theta Z^n_s|^2\ud
s)^{\frac12}\Big]<\infty$ and hence the process $\I e^{\beta s} D_\theta
Y^n_s D_\theta Z^n_s \ud \hat{W}_s$ is well defined.

From (\ref{aux:2:10.jan.2007}) we obtain by taking conditional
$Q_n$-expectations \be |D_\theta Y^n_t|^2 &\leq& C
\E^{Q_n}\Big[\,|D_\theta \xi|^2+\int_0^T |D_\theta Y^n_s
|\,|D_\theta l(s,Y^n_s,Z^n_s)|\ud s\,\big|\F_t\Big]. \ee Next,
integrating in $\theta$, using Fubini's Theorem and Doob's $L^p-$
inequality, we get \be
&&\E^{Q_n}\Big[\Big(\sup_{t\in [0,T]} \int_0^T |D_\theta Y^n_t|^2\ud \theta \Big)^p\Big]\\
&&\leq
C\E^{Q_n}\Big[ \sup_{t\in[0,T]}
\Big(\E^Q\Big[\int_0^T|D_\theta \xi|^2\ud \theta+\int_0^T
\int_0^T |D_\theta Y^n_s |\,|D_\theta l(s,Y^n_s,Z^n_s)|\ud \theta\ud s\,\big|\F_t\Big] \Big)^p\Big]\\
&&\leq C\E^{Q_n}\Big[ \Big(\int_0^T|D_\theta \xi|^2\ud
\theta\Big)^p+ \Big(\int_0^T\int_0^T |D_\theta Y^n_s |\,|D_\theta
l(s,Y^n_s,Z^n_s)|\ud \theta\ud s\Big)^p\Big]. \ee
The last term on the right hand side of the preceding inequality can
be simplified using H\"older's and Young's inequalities with the result \be
&&\int_0^T\int_0^T  |D_\theta Y^n_s |\,|D_\theta l(s,Y^n_s,Z^n_s)|\ud \theta\ud s\\
&&\qquad\leq \int_0^T \Big[\Big(\int_0^T |D_\theta Y^n_s|^2\ud \theta\Big)^{\frac12}\Big(\int_0^T |D_\theta l(s,Y^n_s,Z^n_s)|^2\ud \theta\Big)^{\frac12}\Big]\ud s\\
&&\qquad\leq \sup_{t\in[0,T]} \Big(\int_0^T |D_\theta Y^n_t|^2\ud \theta\Big)^{\frac12}
\int_0^T \Big(\int_0^T |D_\theta l(s,Y^n_s,Z^n_s)|^2\ud \theta\Big)^{\frac12}\ud s \\
&&\qquad\leq
 \frac{1}{C_1}\sup_{t\in[0,T]}\int_0^T |D_\theta Y^n_t |^2\ud \theta
 +C_2\int_0^T\int_0^T |D_\theta l(s,Y^n_s,Z^n_s)|^2\ud \theta\ud s.
\ee
Since for $a,b\geq 0$ we have $(a+b)^p\leq C_3(a^p+b^p)$, by choosing $C_1$ conveniently we obtain
\ben
&&\E^{Q_n}\Big[\Big(\sup_{t\in [0,T]} \int_0^T |D_\theta Y^n_t|^2\ud \theta \Big)^p\Big]\nonumber\\
&&\qquad\leq
C_4\E^{Q_n}\Big[ \Big(\int_0^T|D_\theta \xi|^2\ud \theta\Big)^p
 +\Big(\int_0^T\int_0^T |D_\theta l(s,Y^n_s,Z^n_s)|^2\ud s\ud \theta\Big)^p\Big],\label{aux:3:10.jan.2007}
\een
which provides the desired bound for the part of the Malliavin derivatives of $Y^n$.
Concerning the inequality for the Malliavin derivatives of the $Z^n$ part,
we consult again Equation (\ref{aux:2:10.jan.2007}), from which we derive
\be
\int_0^T \ebs|D_\theta Z^n_s|^2\ud s&\leq& \ebT |D_\theta \xi|^2-2\int_0^T \ebs
D_\theta Y^n_s D_\theta Z_s^n \ud \hat{W}_s
\\
&&
+2\int_0^T \ebs |D_\theta Y^n_s||D_\theta l(s,Y^n_s,Z^n_s)|\ud s
+2\int_0^T M\ebs |D_\theta Y^n_s||D_\theta Z^n_s|\ud s.
\ee
Further estimate
\be
2\int_0^T M\ebs |D_\theta Y^n_s||D_\theta Z^n_s|\ud s&\leq&
 4M^2\int_0^T \ebs |D_\theta Y^n_s|^2\ud s + \frac{1}{2}\int_0^T \ebs |D_\theta Z^n_s|^2\ud s,\\
2\int_0^T \ebs |D_\theta Y^n_s||D_\theta l(s,Y^n_s,Z^n_s)|\ud s&\le&
\int_0^T \ebs |D_\theta Y^n_s|^2\ud s + \int_0^T \ebs |D_\theta l(t,Y^n_s,Z^n_s)|^2\ud s.
\ee
Hence the initial estimate leads to
\be
\frac12\int_0^T \ebs|D_\theta Z^n_s|^2\ud s&\leq&  \ebT|D_\theta \xi|^2-2\int_0^T \ebs
D_\theta Y^n_s D_\theta Z_s^n \ud \hat{W}_s
\\
&&+(1+4M^2)\int_0^T \ebs |D_\theta Y^n_s|^2\ud s+ \int_0^T \ebs |D_\theta l(s,Y^n_s,Z^n_s)|^2\ud s.
\ee
Now for $p>1$ integrate in $\theta$, take $Q_n-$expectations, using Fubini's Theorem as well as a stochastic
version of it to estimate
\ben
&&\E^{Q_n}\Big[\Big(\int_0^T\int_0^T |D_\theta Z^n_s|^2\ud \theta\ud s\Big)^p\Big]\nonumber\\
&&\leq C_5\Big\{\E^{Q_n}\Big[ \Big(\int_0^T |D_\theta \xi|\ud \theta\Big)^p+\Big(\sup_{t\in[0,T]} \int_0^T |D_\theta Y^n_t|^2\ud \theta\Big)^p\nonumber\\
&&+\Big(\int_0^T \int_0^T |D_\theta l(s,Y^n_s,Z^n_s)|^2\ud \theta \ud s\Big)^p+
\Big(\int_0^T \int_0^T D_\theta Y^n_s\,D_\theta Z^n_s\ud \theta \ud \hat{W}_s\Big)^p
\Big]\Big\}.\label{aux:4:10.jan.2007}
\een
We estimate the last term using Burkholder-Davis-Gundy's inequality, which results in
\be
\E^{Q_n}\Big[\Big(\int_0^T \int_0^T D_\theta Y^n_s D_\theta Z^n_s \ud \theta \ud \hat{W}_s\Big)^p\Big]
&\le&
C_6 \E^{Q_n}\Big[\Big(\int_0^T \Big\{\int_0^T |D_\theta Y^n_s||D_\theta Z^n_s|\ud
\theta\Big\}^2 \ud s\Big)^{\frac p 2}\Big].
\ee
Using Cauchy-Schwarz' inequality, we estimate further by
\[
\int_0^T |D_\theta Y^n_s||D_\theta Z^n_s|\ud \theta \leq
\Big(\int_0^T |D_\theta Y^n_s|^2\ud \theta\Big)^{\frac12}
\Big(\int_0^T |D_\theta Z^n_s|^2\ud \theta\Big)^{\frac12}.
\]
Then, with another application of Young's inequality, we obtain \be
&&\E^{Q_n}\Big[\Big(\int_0^T \Big\{\int_0^T |D_\theta Y^n_s||D_\theta Z^n_s|\ud \theta\Big\}^2 \ud s\Big)^{\frac p 2}\Big]\\
&&\qquad\leq
\E^{Q_n}\Big[\Big(\int_0^T \Big[\Big\{\int_0^T |D_\theta Y^n_s|^2\ud \theta\Big\}^{\frac12} \Big\{\int_0^T |D_\theta Z^n_s|^2\ud \theta\Big\}^{\frac12}\Big]^2 \ud s\Big)^{\frac p2}\Big]\\
&&\qquad\leq
\E^{Q_n}\Big[
\Big\{\sup_{t\in[0,T]} \int_0^T |D_\theta Y^n_t|^2\ud \theta \Big\}^{\frac p2}\,
\Big\{\int_0^T  \int_0^T |D_\theta Z^n_s|^2\ud \theta \ud s\Big\}^{\frac p2}\Big]\\
&&\qquad\leq
\frac{1}{C_7}\E^{Q_n}\Big[
\Big\{\sup_{t\in[0,T]} \int_0^T |D_\theta Y^n_t|^2\ud \theta \Big\}^{p}\Big]+
C_7\E^{Q_n}\Big[\Big\{\int_0^T  \int_0^T |D_\theta Z^n_s|^2\ud \theta \ud s\Big\}^{p}\Big].
\ee
Using this last estimate with $C_7$ chosen properly in conjunction with (\ref{aux:3:10.jan.2007})
in (\ref{aux:4:10.jan.2007}),
we obtain
\ben
&&\E^{Q_n}\Big[\Big(\int_0^T\int_0^T |D_\theta Z^n_s|^2\ud \theta\ud s\Big)^p\Big]\nonumber\\
&&\qquad\leq
C_8\E^{Q_n}\Big[ \Big(\int_0^T |D_\theta \xi|\ud \theta\Big)^p
+\Big(\int_0^T \int_0^T |D_\theta l(s,Y^n_s,Z^n_s)|^2\ud \theta \ud s\Big)^p
\Big]\label{aux:5:10.jan.2007}.
\een
Combining inequalities (\ref{aux:3:10.jan.2007}) and (\ref{aux:5:10.jan.2007}) yields the desired estimate and proves the Lemma.
\end{proof}

In the same fashion as in Section \ref{secesti}, we can now
combine the result of the a priori inequality under $Q_n, n\in\N,$
with the inverse H\"older inequality in disguise of Lemma
\ref{bmoeigen} to upgrade the a priori estimates to the following
one. In fact, we observe $|g'_n(z)|\leq 2|z|$ for $z\in\R, n\in\N$.
Moreover, a careful analysis of the demonstration of Lemma 1 of
\cite{06Mor} shows that each $\int Z^n dW$ is also $BMO$ and there
exists a constant $K$ such that \ben \label{uniform:inverse}
\sup_{n\in\N}\lVert \int g'_n(Z^n)
dW\lVert_{BMO_2}\le\sup_{n\in\N}
 \lVert \int Z^n dW\lVert_{BMO_2}+\lVert \int Z
dW\lVert_{BMO_2}=K<\infty. \een  So by Lemma \ref{bmoeigen} there
exists a $1<r$ such that $\cE(\int \alpha g'_n(Z^n) dW)_T$ and
$\cE(\int \alpha Z dW)_T$ are in $L^r(P)$ for all $n\in \N$ with
$r$-norms bounded in $n$. So, again we may apply the argument
based on the third statement of Lemma \ref{bmoeigen}, uniformly in
$n$. This, together with a similar argument applied to the
sequence $(Y^n, Z^n)$ leads to the following a priori estimate.

\begin{lemma}\label{apriori:malliavin:upgrade}
Let $p > 1$ and $r>1$ such that $\mathcal{E}(\int \alpha g'_n(Z^n)
dW)_T \in L^r(P)$ for all $n\in\N$ with a uniform bound. Then
there exists a constant $C > 0$, depending only on $p$, $T$ and
$K$ (from \ref{uniform:inverse}), such that with the conjugate
exponent $q$ of $r$ we have \be &&\E^{P}
\Big[\int_0^T\int_0^T|D_\theta Y^n_t|^{2}\ud \theta \ud
t\Big]^p+\E^{P}
\left[ \int_0^T \int_0^T|D_\theta Z^n_s|^2\ud \theta\ud s  \right]^p\\
&&\quad\leq
    C\Big\{\E^{P}\Big[\Big(\int_0^T|D_\theta \xi|^{2}\ud \theta +
    \int_0^T\int_0^T | \tilde{K}_\theta(\omega,t)|^2\ud \theta\ud t\Big)^{p q^2}\Big]^{\frac{1}{q^2}}\\
&&\quad+
\E^P\Big[\Big(\sup_{t\in[0,T]}\int_0^T |K_\theta(\omega,t)|^2\ud \theta\Big)^{2pq^2}\Big]^{\frac{1}{2q^2}}
\E^P\Big[\Big(|\xi|^{2}+\int_0^T |l(\omega, t, 0,0)|^2\ud s\Big)^{4pq^4}\Big]^{\frac{1}{2q^4}}
    \Big\}<\infty.
\ee
\end{lemma}

\begin{proof}
The proof of the lemma is achieved in three steps.

We start by applying the third statement of Lemma \ref{bmoeigen} as
in the proof of Theorem \ref{apriori.estimate} to the result of
Lemma \ref{lemma.aux:2:11.jan.2007}, from which we obtain with a
constant $C_1$ not depending on $n$ \be &&\E^{P}
\Big[\int_0^T\int_0^T|D_\theta Y^n_t|^{2}\ud \theta \ud
t\Big]^p+\E^{P}
\left[ \int_0^T \int_0^T|D_\theta Z^n_s|^2\ud \theta\ud s  \right]^p\\
&&\qquad\leq
    C_1 \E^{P}\Big[\Big(\int_0^T|D_\theta \xi|^{2}\ud \theta +
    \int_0^T\int_0^T |D_\theta l(s,Y^n_s,Z^n_s)|^2\ud \theta\ud s\Big)^{p q^2}\Big]^{\frac{1}{q^2}}.
\ee

In a second step, we have to estimate the last term of the
preceding equation. From Condition (E2) we obtain with another
universal constant \be
&&\int_0^T\int_0^T |D_\theta l(s,Y^n_s,Z^n_s)|^2\ud \theta\ud t\\
&&\leq
C_2\int_0^T\int_0^T | [\tilde{K}_\theta(\omega,t)|^2+|K_\theta(\omega,t)|^2(|Y^n_t|^2+|Z^n_t|^2)]\ud \theta\ud t\\
&&\leq C_2\Big\{\int_0^T\int_0^T | \tilde{K}_\theta(\omega,t)|^2\ud
\theta\ud t +\sup_{t\in[0,T]}\int_0^T |K_\theta(\omega,t)|^2\ud
\theta \int_0^T |Y^n_s|^2+|Z^n_s|^2\ud s\Big\}. \ee Hence an
application of H\"older's inequality results in \ben
&&\E[\Big(\int_0^T\int_0^T |D_\theta l(s,Y^n_s,Z^n_s)|^2\ud
\theta\ud s\Big)^{pq^2}] \leq
C_3\E[\Big(\int_0^T\int_0^T | \tilde{K}_\theta(\omega,t)|^2\ud \theta\ud s\Big)^{pq^2}]\nonumber\\
&&\qquad+ C_4\E[\Big(\sup_{t\in[0,T]}\int_0^T
|K_\theta(\omega,t)|^2\ud \theta\Big)^{2pq^2}]^{\frac12}
\E[\Big(\int_0^T |Y^n_s|^2+|Z^n_s|^2\ud
s\Big)^{2pq^2}]^{\frac12}.\label{eq:aux:2:24.Jan.2007} \een

In a last step, we need to provide a bound for the $\E[\Big(\int_0^T
|Y^n_s|^2+|Z^n_s|^2\ud s\Big)^{2pq^2}]$ term. For this purpose,
we shall use another application of Theorem
\ref{apriori.estimate}, uniformly in $n$. It requires the
intervention of a different family of measure changes depending on
$n$, which can again be controlled by the BMO property of the
intervening martingales and the third statement of Lemma
\ref{bmoeigen}. In fact, comparing (\ref{linbsde.aux:24.jan.2007})
with (\ref{bsde-lin-001}), we see that the analogue of $H$ has to
be given by $h_n(Z^n)$, where
\[
h_n(z):=\frac{g_n(z)}{z}=\left\{
\begin{array}{cl}
z &,\textrm{if } |z|\leq n\\
\frac{2n|z|-n^2}{z}&,\textrm{if }|z|>n,
\end{array}
\right.
\] which is obviously well defined for all $z$.
In this situation, the stochastic integrals of $h_n(Z^n)$ generate
$BMO$ martingales with uniformly bounded norms. More precisely,
since $g_n(z)\leq z^2$, we have $\sup_{n\in\N}|h_n(z)| \leq |z|$,
$z\in\R.$ A careful analysis of the demonstration of Lemma 1 of
\cite{06Mor} shows that each $Z^n$ is also $BMO$ and there exists
a constant $K$ such that \be \sup_{n\in\N}\lVert \int Z^n
dW\lVert_{BMO_2}+\lVert \int Z dW\lVert_{BMO_2}=K<\infty. \ee Due
to the definition of $h_n$, we may extend (\ref{uniform:inverse})
to \be
 \sup_{n\in\N}\lVert \int h_n(Z^n)
\ud W\lVert_{BMO_2} \leq \sup_{n\in\N} \lVert \int Z^n \ud
W\lVert_{BMO_2}\le K. \ee
By Lemma \ref{bmoeigen} we may assume that $\cE(\int \alpha h_n(Z^n) \ud W)_T$ and $\cE(\int
\alpha Z \ud W)_T$ are in $L^r(P)$ for all $n\in \N$ with $r$-norms
bounded in $n$. 
So, again we may apply the argument based on the third
statement of Lemma \ref{bmoeigen}, uniformly in $n$. We obtain
from Theorem \ref{apriori.estimate} with the settings
$A_t=l(\omega, t, 0,0)$, $\tilde{l}(\omega, t,
Y^n_t,Z^n_t):=l(\omega, t, Y^n_t,Z^n_t)-l(\omega, t, 0,0)$, $\zeta
= \xi$ and $H_t=h_n(Z^n_t)$, for all $\gamma > 1$ the inequality
\ben\label{eq.aux:1:25.jan.2007} \sup_{n\in
\N}\E^P\Big[\Big(\sup_{t\in[0,T]} |Y^n_t|^2 + \int_0^T
|Z^n_s|^2\ud s\Big)^{\gamma}\Big]\leq C_4 \E^P\Big[ \Big(
|\xi|^2+\int_0^T |l(\omega, s, 0,0)|^2\ud s\Big)^{\gamma q^2}
\Big]^{\frac{1}{q^2}}. \een Plugging this inequality into
(\ref{eq:aux:2:24.Jan.2007}) terminates the proof of the Lemma.
\end{proof}

Our main result can now be proved.

\begin{theo}
\label{theo.mall.diff} Assume that driver and terminal condition
satisfy hypotheses (E\ref{para3}), (E\ref{para5}) and
(E\ref{para4}). Then the solution processes $(Y_t,Z_t)$, of
(\ref{linbsde.aux015}) belongs to
$\mathbb{D}^{1,2}\times(\mathbb{D}^{1,2})^d$ and a version of
$(D_\theta Y_t, D_\theta Z_t)$ satisfies for $1\leq i\leq d$\ben
D^i_\theta Y_t &=& 0, \qquad D^i_\theta Z_t = 0,\qquad t\in [0, \theta),\nonumber\\
D^i_\theta Y_t &=& D^i_\theta \xi + \int_t^T \Big[\partial_y
l(s,Y_s,Z_s)D^i_\theta Y_s+\partial_z l(s,Y_s,Z_s)D^i_\theta Z_s
+D^i_\theta l(s,Y_s,Z_s) +2\alpha Z_s D^i_\theta Z_s \Big]\ud s\nonumber\\
&& - \int_t^T D^i_\theta Z_s \ud W_s, \qquad t\in
[\theta,T].\label{mal.dif} \een Moreover, $\{D_t Y_t:\, 0\leq t\leq
T \}$ is a version of $\{Z_t:\, 0\leq t\leq T\}$.

\end{theo}

\begin{proof}
Again, we simplify notation by just considering the case of one
dimensional $Z$ and Wiener process. We first apply Lemma
\ref{apriori:malliavin:upgrade} to obtain a bound of the $L^2$-norms
of the processes $(D_\theta Y^n_t, D_\theta Z^n_t)$, which is
uniform in $n$. Lemma
\ref{apriori:malliavin:upgrade} furthermore allows an appeal to a
weak compactness result to deduce the existence of a pair of
processes $(U_{\theta,t},V_{\theta,t}), 0\le \theta, t\le T,$ and
a subsequence $(n_i)$ such that $(D_\theta Y^{n_i}_t(\omega),
D_\theta Z^{n_i}_t(\omega))$ converges to $(U_{\theta,t}(\omega),V_{\theta,t}(\omega))$ in the weak topology of
the $L^2$ space of random variables with values in $L^2([0,T]\times[0,T])$.

For almost all $t$, Lemma \ref{lemma.from.nualart} implies that $(Y_t,Z_t)$ is Malliavin
differentiable and the equality $(D_\theta Y_t, D_\theta
Z_t)=(U_{\theta,t},V_{\theta,t})$ holds almost everywhere in
$\Omega\times[0,T].$

It remains to use these convergence properties to deduce convergence
term by term in (\ref{aux:1:10.jan.2007}) to (\ref{mal.dif}).

We first show that the stochastic integral terms converge weakly in $L^2(\Omega)$. To this end let $\Psi \in L^2(\Omega)$ be $\cF_T$-measurable. Then there exists a predictable $\psi \in L^2(\Omega \times [0,T])$ with $\Psi = \E(\Psi) + \int_0^T \psi_s d W_s$, and hence 
\be
\lim_{i\to \infty} \E\left[ \Psi \int_0^T D_\theta Z^{n_i}_s dW_s \right] &=& \lim_{i\to \infty} \E\left[ \int_0^T \psi_s D_\theta Z^{n_i}_s d s \right]  = \E\left[ \int_0^T \psi_s D_\theta Z_s ds \right] \\
&=&\E\left[ \Psi \int_0^T D_\theta Z_s dW_s \right],
\ee
which shows that $\int_0^T D_\theta Z^{n_i}_s dW_s$ converges weakly to $\int_0^T D_\theta Z_s dW_s$
in $L^2(\Omega\times[0,T])$.

Next observe that for any bounded $\cF_T$-measurable random variable $B$ we have
\ben\label{alibaba}
&& \E[ B \int_0^T (g'_{n_i}(Z^{n_i} _t)D_\theta Z^{n_i}_t-g'(Z_t)D_\theta Z_t) \ud t] \nonumber\\
&&\qquad= \E[ B \int_0^T (g'_{n_i}(Z^{n_i}_t)-g'(Z_t))D_\theta Z^{n_i}_t \ud t] + \E[ B \int_0^T g'(Z_t)(D_\theta Z^{n_i}_t - D_\theta Z_t) \ud t] 
\een
The first summand on the RHS of Equation (\ref{alibaba}) is bounded by
\be
(\textrm{esssup}\ |B|) \ \sup_i \left(\E\int_0^T D_\theta (Z^{n_i}_t)^2 \ud t\right)^\frac12 \left(\E[\int_0^T (g'_{n_i}(Z^{n_i}_t)-g'(Z_t))^2 \ud t]\right)^\frac12,
\ee
which converges to $0$ as $i \to \infty$. The second summand on the RHS of Equation (\ref{alibaba}) converges also to $0$ since $B g'(Z_t)\in L^2(\Omega\times[0,T])$ and $D_\theta Z^{n_i}_t$ converges weakly to $D_\theta Z_t$. Since $B$ was arbitrary we have shown that $\int_0^T g'_{n_i}(Z^{n_i} _t)D_\theta Z^{n_i}_t\ud t$ converges to $\int_0^Tg'(Z_t)D_\theta Z_t\ud t$ in the weak topology of $L^1(\Omega \times [0,T])$. 
%

Finally we come to the various derivative terms of $l$. The boundedness of
the partial derivatives of $l$ as well as Condition (E\ref{para5}) and Inequality (\ref{eq.aux:1:25.jan.2007}) imply
\be
& &\sup_n \E \int_0^T \int_0^T [D_\theta (l(s, Y^n_s,Z^n_s))]^2 \ud s\, \ud \theta \\
&\le& C_1 \sup_n \E \int_0^T \int_0^T |\partial_y l(Y^n_s,Z^n_s)D_\theta
Y^n_s|^2+ |\partial_z l(Y^n_s,Z^n_s)D_\theta Z^n|^2 + [(D_\theta l)(s, Y^n_s,Z^n_s)]^2 \ud s\, \ud \theta \\
&\le& C_2 \Big\{ \sup_n \E \int_0^T \int_0^T |D_\theta Y^n_s|^2+ |D_\theta Z^n_s|^2 + |\tilde K_\theta(s)|^2 \ud s\, \ud \theta \\
& & \qquad + \sup_n \E \left[\int_0^T \sup_{s\in[0,T]}|K_\theta(s)|^2\ud \theta  \int_0^T (|Y^n_s|+|Z^n_s|)^2 \ud s \right] \Big\} < \infty.
\ee
Thus, by Lemma \ref{lemma.from.nualart}, for almost all $s\in[0,T]$, $l(s,Y_s, Z_s)$ belongs to $\mathbb D^{1,2}$, and 
$D_\theta(l(s,Y^n_s, Z^n_s))$ converges to $D_\theta(l(s,Y_s, Z_s))$ weakly in $L^2(\Omega\times [0,T])$. Since the partial derivatives of $l$ are continuous we have $D_\theta(l(s,Y_s, Z_s)) = \partial_y l(Y_s,Z_s)D_\theta
Y_s+ \partial_z l(Y_s,Z_s) D_\theta Z_s + (D_\theta l)(s, Y_s,Z_s)$. 
%
%
\end{proof}

We next assume more structural properties for the terminal
variable. More precisely, we will turn to the framework of
forward-backward systems. 
Given a $d-$dimensional Brownian motion $W$ and an $x\in\R^n$, we denote by
$X_t=(X_t^1,\cdots,X_t^n)$ the forward part given by \ben\label{sde.mall}
X_t=x+\int_0^t b(s,X_s)ds+\int_0^t [\sigma(s,X_s)]^*\ud W_s .\een  The
coefficients are supposed to satisfy
\begin{enumerate}
\item[(P1)]
 $b,\sigma_i:[0,T]\times\R^n\to\R^n$,  are uniformly Lipschitz; $b(\cdot,0)$ and $\sigma_i(\cdot,0)$ are bounded for $1\leq i\leq d$;
 $\sigma(t,X_t)$ is a $d\times n$ matrix and $[\sigma(t,X_t)]^*$ represents its transpose.
\end{enumerate}

For the backward part we consider \ben\label{eq.aux:30.jan.2007}
Y_t = g(X_T)- \int_t^T Z_s dW_s + \int_t^T
[l(s,X_s,Y_s,Z_s)+\alpha|Z_s|^2] ds, \qquad t\in[0,T] \een where
the driver and the terminal conditions are supposed to satisfy the
following assumptions
\begin{enumerate}
\item[(P2)] $l:[0,T]\times\R^n\times\R\times\R^d\to \R$ is globally
Lipschitz and continuously differentiable in $(x,y,z)$;
$l(\cdot,0,0,0)$ and $\partial_x l(t,x,y,z)$ are bounded by a constant $M$.

\item[(P3)] $g:\R^n\to\R$ is a bounded differentiable function with
bounded first partial derivatives.
\end{enumerate}

Before stating the main theorem we recall that the processes $X,Y$
and $Z$ all depend on the variable $x$. For ease of notation we
omit the corresponding superscripts. In this setting our main
result is the following.

\begin{theo}\label{theo.mall.FBSDE}
Suppose that the coefficients of the SDE (\ref{sde.mall}) and the
driver and terminal condition of the BSDE
(\ref{eq.aux:30.jan.2007}) satisfy conditions (P1), (P2) and (P3).
Then the solution processes $(X,Y,Z)$
possess the following properties.
\begin{itemize}
\item
For any $0\leq t\leq T$, $x\in\R$, $(Y_t,Z_t)\in
\mathbb{D}^{1,2}\times\big(\mathbb{D}^{1,2}\big)^d$,  and a version of
$\{(D^i_\theta Y_t,D^i_\theta Z_t);0\leq \theta,t\leq T\}$ satisfies for $1\leq i\leq d$
\ben
D^i_\theta Y_t &=& 0, \qquad D^i_\theta Z_t = 0,\qquad t\in [0, \theta),\nonumber\\
D^i_\theta Y_t &=& \partial_x g(X_T) D^i_\theta X_T + \int_t^T \Big[\partial_x l(s,Y_s,Z_s)D^i_\theta X_s+\partial_y l(Y_s,Z_s)D^i_\theta Y_s\nonumber\\
&&+\partial_z l(Y_s,Z_s)D^i_\theta Z_s +2\alpha Z_s D^i_\theta Z_s \Big]\ud s
 - \int_t^T \langle D^i_\theta Z_s, \ud W_s\rangle, \qquad t\in
[\theta,T].\label{eq.aux:40} \een Moreover, $\{D_t Y_t; 0\le t
\leq T\}$ defined by the above equation is a version of $\{Z_t;
0\le t\leq T\}$. \item The following set of equations holds for
any $0\leq\theta\leq t\leq T$ and $x\in\R^n$, $P$-almost surely,
\be
D_\theta X_t &=& \partial_x X_t (\partial_x X_\theta)^{-1}\sigma(\theta,X_\theta)\\
D_\theta Y_t &=& \partial_x Y_t (\partial_x X_\theta)^{-1}\sigma(\theta,X_\theta)\\
Z_t&=&\partial_x Y_t (\partial_x X_t)^{-1}\sigma(s,X_t); \ee
and $D_\theta Z_t =\partial_x Z_t (\partial_x X_\theta)^{-1}\sigma(\theta,X_\theta)$ for almost all $(\om, t)$. 
\end{itemize}
\end{theo}

\begin{proof}
Theorem 2.2.1 of \cite{nualart} assures existence, uniqueness and
Malliavin differentiability of solutions of  SDE (\ref{sde.mall})
under Hypothesis (P1). Moreover the solution processes satisfy
$X_t\in(\mathbb{D}^{1,\infty})^n$ for any $t\in[0,T]$ and $1\leq i\leq
d$ and the following equation holds: \be
D^i_\theta X_t &=& 0,\phantom{9877ncfgryusfngncerif8ync97st59nsc8e75yt9ns8chgfhgscdu} t\in [0, \theta),\\
D^i_\theta X_t &=& \sigma(\theta, X_\theta)  + \int_\theta^t
\partial_x b(s,X_s)D^i_\theta X_s\ud s +\int_\theta^t
\partial_x \sigma (s, X_s)D^i_\theta X_s\ud W_s, \qquad t\in
[\theta,T], \ee
(see f.ex. Theorem 2.2.1 of \cite{nualart}).

Let us next check the validity of hypotheses (E1)-(E3) for the
driver of our BSDE, for simplicity in the one dimensional case.
From Condition (P3) it follows that the function $g$ and its
derivative are bounded. In combination with the fact that $X_t\in
\mathbb{D}^{1,\infty}$ this implies that
$g(X_T)\in\mathbb{D}^{1,\infty}$, i.e. in the setting of Theorem
\ref{theo.mall.diff} Condition (E3) is verified. From Condition
(P2), we have $|l(t,X_t,0,0)|\leq M(1+|X_t|)$ for $t\in[0,T].$ The
fact that $X_t\in\mathbb{D}^{1,\infty}$ then entails
$\E[\Big(\int_0^T |l(t,X_t,0,0)|^2\ud s\Big)^p]<\infty$ for all
$p\geq 1$. Hence (E1) is satisfied. Condition (P2) includes the
statement that $\partial_x l(\cdot,x,\cdot,\cdot)$ is bounded.
Therefore we have
\[
|D_\theta l(t,X_t,y,z)|=|\partial_x l(t,X_t,y,z)D_\theta X_t|\leq
M |D_\theta X_t|
\]
with some constant $M$. Using the fact that
$X_t\in\mathbb{D}^{1,\infty}$ we obtain $\E[\Big(\int_0^T \int_0^T
|M\,D_\theta X_t|^2\ud \theta \ud t \Big)^p]<\infty$ for $p>1$,
which means Condition (E2) also holds.

With conditions (E1), (E2) and (E3) verified we can apply Theorem
\ref{theo.mall.diff}, which implies the Malliavin
differentiability of $(Y_s,Z_s)$ and proves the first block of
results.

For the second part of the theorem, the representation formula of
$D X$ is standard (see f.ex.\ Chapter 2.3 in \cite{nualart}). The representation of $Z$ by the trace of $D Y$
being granted, we only have to prove the representation formulas
for $D Y$ and $D Z$. For this purpose, we apply It\^o's formula to
$\partial_x Y_t (\partial_x
X_\theta)^{-1}\sigma(\theta,X_\theta)$, then use
(\ref{pararealdiff}) to represent the $\partial_x Y_t$ term. We
further use the representation of $DX$ to account for the terminal
condition. This way we obtain (\ref{eq.aux:40}) with $D_\theta
Y_t=\partial_x Y_t (\partial_x
X_\theta)^{-1}\sigma(\theta,X_\theta)$ and $D_\theta
Z_t=\partial_x Z_t (\partial_x
X_\theta)^{-1}\sigma(\theta,X_\theta)$. The representation follows
from uniqueness of solutions for the BSDE.
\end{proof}
{\bf Example:}\\We finally study a specific setting of Theorem
\ref{theo.mall.FBSDE}. We assume that $\xi$ and $D \xi $ are bounded
by $M$. Assume further that the driver does not depend on $y$. Then,
choosing $\theta=t$, representation (\ref{eq.aux:40}) can be
simplified to \be Z_t &=& D_t \xi + \int_t^T \Big[\partial_z
l(Z_s)D_t Z_s +2\alpha Z_s D_t Z_s \Big]\ud s - \int_t^T D_t Z_s \ud
W_s. \ee Since $\int Z dW\in BMO$ and if we further assume $\int
\partial_z l(Z) dW\in BMO$ we may change the measure to
$Q=\cE\Big(\int \partial_z l(Z)+2\alpha Z dW_T\Big)\cdot P$. Hence
we obtain, by applying conditional expectations \be Z_t = \E^Q[D_t
\xi|\F_t]\leq M \E^Q[1|\F_t]\leq M. \ee This means $Z\in L^\infty$.
This way we recover the Malliavin differentiability results of
\cite{06HM} from our main result.
\newline
 \textbf{ Remarks:}

1. The methods of proof of this Section, building upon a truncated
sequence of Lipschitz BSDEs, could also be used in the treatment of
the differentiability problem in Section \ref{sec:proof:diff}.
This sequence would allow the use of the results in \cite{97KPQ},
which, combined with the a priori estimates of sections
\ref{secesti} and \ref{secpriori} would imply differentiability.

2. Our main results allow less restrictive hypotheses. For example
in Section \ref{secesti}, we assume for our a priori estimations that
$\zeta\in L^p$ for all $p\geq1$. An analysis of the proof clearly
reveals that to obtain estimates in $\mathcal{R}^p$ or
$\mathbb{L}^p$ we only need that $\zeta\in L^p$ for all
$p\in(2,2pq^2]$. We chose to write $\zeta\in L^p$ for all $p\geq
1$ not to produce an overload of technicalities in a technically
already rather complex text.

{\bf Acknowledgement:} We thank an anonymous referee for a very careful reading and for many helpful remarks.
This research was supported by the DFG Research Center MATHEON "Mathematics for Key Technologies" (FZT86) in Berlin. 


\end{document}